\documentclass[11pt]{article}
\input epsf
\usepackage{enumerate, fullpage}
\usepackage{array}
\usepackage{graphicx}
\usepackage{multirow}
\usepackage{color}
\usepackage{amsfonts}

\def\b1{\mathbf{1}}

\def\cO{{\cal O}}

\newcommand{\field}[1]{\mathbb{#1}}

\newcommand{\diag}{\mathop{\bf diag}}
\newcommand{\beq}{\begin{equation}}
\newcommand{\eeq}{\end{equation}}

\newcommand{\beqnr}{\begin{eqnarray}}
\newcommand{\eeqnr}{\end{eqnarray}}

\newcommand{\benum}{\begin{enumerate}}
\newcommand{\eenum}{\end{enumerate}}
\newcommand{\argmax}{\mathop{\rm argmax}}

\newcommand{\QED}{\rule{.1in}{.1in}}
\newcommand{\cA}{{\cal A}}

\newcommand{\cK}{{\cal K}}

\newcommand{\cG}{{\cal G}}

\newcommand{\cD}{{\cal D}}

\newcommand{\Observe}{{\bf Observe\ }}

\newcommand{\For}{{\bf For\ }}
\newcommand{\Do}{{\bf do\ }}

\newcommand{\D}[1]{{\mathbb#1}}
\newcommand{\Rn}{{\D{R}^n}}

\newtheorem{DE}{Definition}[section]

\newtheorem{NO}[DE]{Notation}

\newtheorem{LE}[DE]{Lemma}

\newtheorem{RE}[DE]{Remark}
\newtheorem{THM}[DE]{Theorem}
\newtheorem{TEMP}[DE]{Template}
\newtheorem{RULE}[DE]{Rule}

\newtheorem{PRO}[DE]{Procedure}

\newcommand{\qed}{\mbox{}\hspace*{\fill}\nolinebreak\mbox{$\rule{0.7em}{0.7em}$}}

\advance \topmargin by -\headheight
\advance \topmargin by -\headsep

\evensidemargin \oddsidemargin

\begin{document}

\begin{center}

\Large {\bf Optimal adaptive control of cascading power grid failures}\footnote{partially funded by grant DE-SC000267}
\vskip 0.1cm
{\Large Version 2010-Dec-10}
\vskip 0.5cm

\Large

Daniel Bienstock\\
Columbia University \\
New York 
\end{center}

\section{Introduction}\label{intro}
Power grids have long been a source of interesting optimization problems.
Perhaps best known among the optimization community are the unit commitment
problems and related
generator dispatching tasks, see \cite{hobbs}.  However, recent
blackout events have renewed interest on problems related to grid 
vulnerabilities.

A difficult problem that has been widely studied, the $N-K$ problem, concerns
the detection of small cardinality sets of lines or buses whose simultaneous
outage could develop into a significant failure event; see \cite{bie2}, \cite{pinar} and references therein.  This is a hard combinatorial problem which, 
unlike the typical formulations for the unit commitment problem, includes
a detailed model of flows in the grid.  A different set of algorithmic 
questions concern how to react to protect a grid when a significant event has taken place.  This is
the outlook that we take in this paper.

In this context, the central modeling ingredient
is that 
power grids display {\em cascading} behavior.  A cascade is
the process by which components of the grid (especially, power lines) sequentially
become inoperative.  In a catastrophic cascade this process accelerates (it 'snowballs') until
the grid collapses.  A control action must take this multi-step behavior 
into account because a myopic action taken at the start of the process so as to 
immediately arrest the cascade may prove far from optimal.

The computation of an optimal control can be formulated as
a multi-stage mixed-integer programming problem; 
ideally this should be a stochastic or robust formulation.  Furthermore the
formulation will need to include an explicit model of the power flows.  When
dealing with large-scale (or even medium-scale) grids it is likely that such
a formulation will prove extremely intractable.  In addition, it is likely
that the prescribed control will call for counter-intuitive and possibly
impractical actions.  See Section \ref{algorithms} and Appendix \ref{traditional}

In this paper, building on prior models for cascades, we consider an affine, 
adaptive, distributive
control algorithm that is computed at the start of the cascade and deployed
during the cascade.  The control sheds demand as a function of observations
of the state of the grid, with the objective of terminating the cascade with
a minimum amount of demand lost.  The optimization problem handled at the
start of the cascade computes the coefficients in the affine control (one
set of coefficients per demand bus).  The discussion of this approach
starts in Section \ref{adaptive1}; Section \ref{firstset} describes an
initial set of experiments with a simple form of affine control.

Most of our algorithms are first-order methods that compute
local optima (see Section \ref{firstorder}); 
however in a special case which is nevertheless 
of interest we obtain
an exact algorithm that runs in polynomial-time (Section \ref{scalingproblem}). Algorithms that account for stochastics are discussed in Section \ref{stochasticoptimization}.  We present numerical experiments 
with parallel implementations of our algorithms, using as data
a snapshot of the U.S. Eastern Interconnect, with 
approximately $15000$ buses and $23000$ lines.

\section{Notation}\label{notation}

In the {\em linearized} approximation to the power flow problem, we are
given a directed graph $G$ with $n$ buses and $m$ lines (denoted, respectively,
``nodes'' and ``arcs'' in traditional graph theoretic language).  In addition
\begin{itemize}
\item A line $j$ is oriented from its {\em tail}, $t(j)$ to its {\em head}, $h(j)$. The orientation of any line is arbitrary and is simply used for notational
convenience.  The set of lines is denoted by $\cA$.
\item For each line $j$ we are given two positive quantities: its {\em flow limit} $u_{j}$ and its {\em reactance} $x_{j}$.
\item We are given a {\em supply-demand} vector $\beta \in \Rn$ with the following interpretation.  For a bus $i$, if $\beta_i > 0$ then $i$ is a {\em generator} (a source bus) while if
$\beta_i < 0$ then $i$ is a {\em load} (a demand bus) and in that case 
$-\beta_i$ is the {\em demand} at $i$.  The condition $\sum_i \beta_i = 0$ is
assumed to hold.  For a generator bus $i$, we indicate by the constant $\tilde s_i$ the maximum supply of $i$. We denote by $\cG$ denote the set of generators and by $\cD$ the
set of demand buses. 
\end{itemize}
The linearized power flow problem specifies a variable $f_{ij}$ associated 
with each line $j$ and a variable $\phi_i$ associated with each bus 
$i$.  Denoting, for each bus $i$, 
the set of lines oriented out of (into) $i$ by $\delta^{+}(i)$ (resp., $\delta^{-}(i)$), the power flow problem consists in finding a solution to the 
following system of equations:  
\begin{eqnarray}
\sum_{j \in \delta^{+}(i)} f_{j} - \sum_{(j)\in \delta^{-}(i)} f_{j}  & = & \beta_i  \ \ \forall \ \mbox{bus} \ i, \label{balance}\\
\phi_{t(j)} - \phi_{h(j)} - x_j f_{j} & = & 0  \ \ \ \forall \ \mbox{line} \ j . \label{ohm-eq}
\end{eqnarray}
These equations can simply be abbreviated as
$$ N f \ = \ \beta, \ \ \ N^T \phi - X f \ = \ 0,$$
\noindent 
where $N$ denotes the bus-arc incidence matrix of $G$, and 
$X = \diag\{x_{ij}\}$.  
\begin{RE}\label{uniqueness}
It can easily be shown that system (\ref{balance})-(\ref{ohm-eq}) is feasible if and only if $\sum_{i \in K} \beta_i = 0$ for each
component (``island'') $K$ of $G$, and in that case the solution is unique
in the $f$ variables.  
\end{RE}
We stress that the orientation of the lines is arbitrary,
 consequently for a line $j$ we might have $f_{j} < 0$, indicating that
power flows in the reverse direction.  As a final point we note that the flow limits $u_{ij}$ 
do not appear in (\ref{balance})-(\ref{ohm-eq}); consequently, it is possible
for the (unique) solution $f$ to exceed the flow limits, i.e. it is 
possible that $|f_{j}| > u_{j}$ for certain lines $j$.
\section{Cascades}\label{cascades}

In this section we introduce our model of cascading failures of power grids, which draws from
the models in \cite{CarrerasCH02}, \cite{Carreras03}, 
\cite{CarrerasCH04}.  A cascade starts with an initial event, for example the
removal of a number of lines, that alter and compromise the power grid. 
Informally, this
is followed by a sequence of additional line outage events which are interspersed
with simple control mechanisms as the grid adjusts to decreased resources, 
for example decreased generator capacity.  

As a starting point, we have the following template:
\begin{center}
  \fbox{
    \begin{minipage}{0.9\linewidth}
      \hspace*{.9in} \begin{TEMP}{GENERIC CASCADE TEMPLATE}\label{gencasc} \end{TEMP}
      {\bf Input}: a power grid with graph $G$ (post-initiating event). Set $G^1 = G$.\\
      \For $r = 1, 2, \ldots$ \Do \\
        \hspace*{.5in}(comment: round $r$ of the cascade)\\
        \hspace*{.5in}{\bf 1.} Set $f^r = $ vector of power flows in $G^r$.\\
        \hspace*{.5in}{\bf 2.} Set $\cO^{r} = $ set of lines of $G^r$ that become outaged in round $r$. \\
        \hspace*{.5in}{\bf 3.} Set $G^{r+1} = G^{r} - \cO^r$.  Adjust demands and supplies in $G^{r}$.\\
    \end{minipage}
  }
\end{center}
In this template, graph $G^1$ represents the grid we have after the initial
event that causes the cascade.  To make the template complete, we must 
specify a mechanism for determining the set $\cO^r$ in Step 2, and how to
adjust demands and supplies in Step 3.  We tackle the second issue first.

The supply/demand adjustment in Step 3 handles cases where 
the removal of $\cO^{r}$ from $G^{r}$ creates new
components in $G^{r+1}$.  Any  imbalance between supply and demand 
in a component of $G^{r+1}$ must be corrected:  
if, in a component of $G^{r+1}$, total supply exceeds total demand, then
supply must be reduced, and viceversa.  Such actions constitute a form of
control.  We will discuss this item in more detail in Section \ref{adaptive1}.

Now we turn to step Step 2.  We will say that line $j$ is
{\em overloaded} in round $r$ if $|f^r_{j}| > u_{j}$.  
Practical experience shows that overloaded lines are 
likely to be come outaged. Two 
ways to implement this observation are:
\begin{itemize}
\item [(F.1)] Simple deterministic rule: $j \in \cO^r$ if
$|f^r_{j}|/u_{j} > 1$ (alternatively, if $|f^r_{j}|/u_{j} \ge 1$).
\item[(F.2)] Stochastic rule: for each line $j$ there is a function
$\gamma_{j}$ such that if $|f^r_{j}|/u_{j} > 1$ then
$j \in \cO^r$ with probability  $\gamma_{j} (|f^r_{j}|/u_{j})$.  In
\cite{Carreras03}, the function $\gamma_{j}$ takes a fixed value $p$.
\end{itemize}
The two alternate forms of (F.1) are not strictly equivalent but we
assume that from a practical perspective they are essentially identical. 
In any case, the use of (F.1) may not be desirable because it represents a
strict numerical criterion that may be difficult to implement using
standard numerical algorithms -- a possibility is to use for
$u_{j}$ a slightly larger value than the true flow limit.  The constant
$\gamma_{j}$ version of rule (F.2) can be criticized in that we would expect that
higher overloads cause outages with higher probabilities; that is to say,
we would expect that $\gamma_{j}(t) \rightarrow 1$ as $t \rightarrow +\infty$.
The problem of choosing, and calibrating such functions $\gamma_{j}$ is 
significant.

Both rules can additionally be criticized on the grounds that they are 
memory-free;  from a realistic perspective 
a line that was highly overloaded
in round $r - 1$ should be more likely to be outaged in round $r$ than one
that was not.  To address this issue, we assume that for any line $j$ 
we are given a parameter $0 \le \alpha_{j} \le 1$ and define quantities
$\tilde f^{r}_{j}$ by
\begin{eqnarray}
\tilde f^{r}_{j} & = & \alpha_{j} |f^{r}_{j}| \, + \, (1 - \alpha_{j}) \tilde f^{r-1}_{j}, \label{memory-1}
\end{eqnarray}
with $\tilde f^{0}_{j}$ set to the absolute value of the flow on $(j)$
prior to the incident that initiates the cascade.  The $\tilde f^{r}_{j}$
quantities are then used instead of $|f^r_{j}|$ to obtain memory-dependent versions of rules (F.1) and (F.2).  A variation of (\ref{memory-1}) is
\begin{eqnarray}
\tilde f^{r}_{j} & = & \alpha_{j} |f^{r}_{j}| \, + \, (1 - \alpha_{j}) |f^{r-1}_{j}|. \label{memory-2}
\end{eqnarray}
The choice of the parameters
$\lambda_{j}$ depends on the time scale of the cascade, but for robustness
purposes the $\lambda_{j}$ should be treated as noisy.

The memory-dependent versions of rules (F.1) and (F.2) can still give rise to non-smooth behavior and ill-conditioning: for example, solving the power 
flow equations with different solvers can give rise to different cascades.  In
order to lessen this difficulty, we introduce an additional detail in
choosing if a line becomes outaged.  
\begin{center}
  \fbox{
    \begin{minipage}{0.9\linewidth}
      \hspace*{.9in} \begin{RULE}{STOCHASTIC LINE OUTAGE}\label{stochout} \end{RULE}
\vskip 2pt
      {\bf Parameters}: $0 \le \epsilon_r \le 1$ for each round $r$.\\
      {\bf Notation}: refer to Template \ref{gencasc} and equation (\ref{memory-1}).\\
      {\bf Application}: For a line $j$ in $G^r$:
\begin{eqnarray}
\mbox{if} && u_{j} < \tilde f^{r}_{j}, \ \ \mbox{then} \ \ j \in \cO^r, \label{actual1} \\
\mbox{if} && (1 - \epsilon_r)u_{j} < \tilde f^{r}_{j} \le u_{j}, \nonumber \\
&&  \mbox{then} \ \  j \in \cO^r \ \ \mbox{with probability $\frac{1}{2}$}, \label{flipping}\\
\mbox{if} && \tilde f^{r}_{j} \le (1 - \epsilon_r)u_{j}, \ \ \ \mbox{then} \ \ \ j \notin \cO^r. \label{actual3}
\end{eqnarray}
\vskip 2pt
    \end{minipage}
  }
\end{center}
\noindent The random choice in (\ref{flipping}) is an indirect way to incorporate some of the (poorly defined) ``noise'' mentioned above;
additionally, from a mathematical perspective, 
it serves to smooth the cascade process.  
Typically we would have $\epsilon_1 \le \epsilon_2 \le \ldots $,
indicating increasing uncertainty as the cascade progresses. If $\epsilon_r = 0$ for all $r$ we obtain the pure deterministic rule.

Rule \ref{stochout}, and extensions, will be used later in our numerical experiments.  

\section{Control Algorithms}\label{algorithms}
We consider control algorithms designed to stop the cascade after a fixed 
number of rounds with a maximum amount of total demand feasibly satisfied.  In
developing such algorithms we assume that the cascade is initially
slow-paced so that significant computation is possible at time zero (immediately after the initiating event).  While this may not be true for all cascades, it
was true in the case of the 2003 cascade in the Northeast U.S. and Canada
\cite{usc} with (arguably) on the order of one hour elapsing between subsequent
outages at the start. Thus, we assume that the algorithm is computed at time zero, with significant
information available as to the state of the grid; the algorithm will be
applied as the cascade progresses, with no further computation. 
We assume, as a control requirement, 
that there is a final round $R$ in the cascade
at the end of which  no lines can be overloaded.  

The computation of an optimal schedule for demand shedding can be stated as mixed-integer 
optimization problem; see Appendix \ref{traditional}.  However, it is not
clear that such an approach is either computationally feasible or even
desirable (see the discussion in the Appendix).  In this paper 
we will take a different outlook.

\subsection{Adaptive control}\label{adaptive1}
We focus on robust control algorithms that take as input
limited, real-time observations on the state of the grid 
and which prescribe simple control
actions such as distributed load shedding (loss of demand).  Our generic
cascade template (\ref{gencasc}) is modified as follows:

\begin{center}
  \fbox{
    \begin{minipage}{0.9\linewidth}
      \hspace*{.9in} \begin{TEMP}{CASCADE CONTROL}\label{contcasc} \end{TEMP}
\vskip 5pt
      {\bf Input}: a power grid with graph $G$. Set $G^1 = G$.\\

      {\bf Step 0.} {\bf Compute} control algorithm.\\

      \For $r = 1, 2, \ldots, R-1$,  \Do \\
        \hspace*{.2in}(comment: controlled round $r$ of the cascade)\\
        \hspace*{.2in}{\bf 1.} Set $f^r = $ vector of power flows in $G^r$.\\
        \hspace*{.2in}{\bf 2.} \Observe state of grid (from state estimation).\\
        \hspace*{.2in}{\bf 3.} {\bf Apply} control.\\
        \hspace*{.2in}{\bf 4.} Set $g^r = $ vector of resulting power flows in $G^r$.\\
        \hspace*{.2in}{\bf 5.} Set $\cO^{r} = $ set of lines of $G^r$ that become outaged in round $r$. \\
        \hspace*{.2in}{\bf 6.} Set $G^{r+1} = G^{r} - \cO^r$.  Adjust loads and generation in $G^{r}$.\\
\vskip 1pt
      {\bf Termination} (round $R$). If any island of $G^R$ has line overloads, proportionally shed demand in that island until all line overloads are eliminated.\footnote{The criterion
of ``stability'' inherent in the termination step may obviously be incomplete 
when using a more complete model of power flows than the linearized model.}
\vskip 6pt
    \end{minipage}
  }
\end{center}

In this template, steps 0, 2 and 3 are the only ones requiring controller 
actions.  The remaining steps are due to the physics of the grid or 
underlying low-level automatic control steps.  As discussed above, in this
paper we assume that the cascade allows enough time for significant
computation to take place in step 0.  One could conceive of a variant of
the template where step 0 is carried out in advance of an initiating event,
thus obtaining a more general form of control.  However, proceeding as in
the template resolves an exponential number of potential outcomes, 
likely obtaining a simpler and faster step 0 and a more effective control
algorithm.

Having applied the control in a given round, a given (pre-existing) 
component 
of $G^r$ is likely to experience a supply/demand imbalance.
This condition must be removed, which we assume can be undertaken 
through a low-level control mechanism.  This is not a trivial assumption
and if the imbalance is large the rebalancing may be deemed impossible with
existing technology, 
resulting in the loss of all demand in that component. It is 
straightforward to implement such a ``maximum imbalance'' feature
in the above template; however for the
sake of simplicity we will assume that we can always rebalance supply and
demand by proportionally decreasing the output of each generator in a 
given component (again, other rebalancing mechanisms are possible, giving rise to alternative
versions of the template).  Having effected the rebalancing, a new set
of power flows will be instantiated: this is vector $g^r$ in Step 4. Steps 5
and 6 now follow as in the generic cascade template.  In Step 6, some
of the (new) components of $G^{r+1}$ may have an excess of supply over demand
or the other way around, and  again we make the assumption that this 
excess is removed through a proportional scaling mechanism.

To make the above template complete, we need to describe the type of control
we have in mind, including what type of data observations it requires and 
what kind of control actions it specifies.  In terms of the last item, 
many possibilities exist (including modifying the structure of the grid
by e.g. shutting down lines, in which case the terminology in Steps 4-6 is
not strictly correct) but in this paper we will focus on one type of action
which is feasible in practice: ``load shedding'' or the controlled loss of
demand.

In the linearized (DC) power flow models we have variables of just two types: 
power flows and phase angles.  In this paper we concentrate on control 
algorithms that observe real-time quantities related to flows.  Two such
quantities are:
\begin{itemize}
\item [(a)] The maximum overload: $\max_{j \in G^r} \{ |f^r_{j}|/u_{j} \}$.  
\item [(b)] The maximum relative flow variability: $\max_{j \in G^{r-1}} \{ |f^r_{j} - f^{r-1}_{j}|/|f^{r-1}_{j}| \}$, where we assume the maximum is taken
over the lines with $f^{r-1}_{j} \neq 0$.
\end{itemize}
From a practical standpoint, a relevant issue is how accurately (a) or (b) 
can be observed in real time, and whether such measurements can be disseminated
to all buses of the grid.  We assume that in early rounds of a cascade this
is not a fundamentally difficult technological problem.  Nevertheless, for
a given integer $\delta > 0$ we define the {\em radius-}$\delta$ version
of either (a) or (b), in which each bus of the grid is expected to perform
the measurement over all links within $\delta$ hops in the current topology.  Note
that even if $\delta$ is large we are still constraining the measurement performed
by a bus $v$ to take place in the same component as $v$.
We will discuss this topic in greater depth below. 

Putting aside this issue, we propose an affine control policy where at
each round $r < R$, each demand bus $v$ independently adjusts its demand
by making use of a (precomputed) triple $(c^r_v, b^r_v, s^r_v)$ of parameters.
\begin{center}
  \fbox{
    \begin{minipage}{0.9\linewidth}
      \hspace*{.9in} \begin{PRO}{AFFINE CONTROL}\label{affinecontcasc} \end{PRO}
      {\bf Input}: a power grid with graph $G$ (post-initiating event). Set $G^1 = G$.\\
      {\bf 0.} {\bf Compute} triples $(c^r_v, b^r_v, s^r_v)$ for each $r < R$ and $v$.\\
      \For $r = 1, 2, \ldots , R-1$, \Do \\
        \hspace*{.5in}(comment: controlled round $r$ of the cascade)\\
        \hspace*{.5in}{\bf 1.} Set $f^r = $ vector of power flows in $G^r$, 
        and $d^r_v = $ the demand of any bus $v$. \\
        \hspace*{.5in}{\bf 2.} For any demand bus $v$, let $\kappa^r_v$ be its data observation.\\
        \hspace*{.65in} {\bf Apply control:} if $\kappa^r_v > c^r_v$, reset the demand of $v$ to $$\min\{1, \, [b^r_v + s^r_v (c^r - \kappa^r_v)]^+ \}~d^r_v.$$ 
        \hspace*{.5in}{\bf 3.} Adjust generator outputs in each component of
$G^r$ so as to match demand.\\
        \hspace*{.5in}{\bf 4.} Set $\cO^{r} = $ set of lines of $G^r$ that become outaged as a result of the flows\\ \hspace*{.65in} instantiated in Step 4. \\
        \hspace*{.5in}{\bf 5.} Set $G^{r+1} = G^{r} - \cO^r$.  Adjust demands and supplies in $G^{r}$.\\

{\bf Round R.} For any component $K$ of $G^R$, set
$\Psi^R_{K} \, \doteq \, \min\left\{ 1 \, , \, \max_{j \in K} \{ |f^R_{j}|/u_{j} \} \right\}$.\\  
\hspace*{.8in}If $\Psi^R_{K} > 1$, then any bus $v$ of $K$ resets its demand to $d^R_v /\Psi^R_{K}$.

    \end{minipage}
  }
\end{center}
In Procedure \ref{affinecontcasc}, round $R$ handles cascade termination.  Under the linearized power flow
the rescaling guarantees that no line overloads will exist.
Note that Step 0 requires the computation of $3(R-1)D$ parameters, where
$D$ is the number of demand buses.  Special cases of the control are:
\begin{itemize}
\item [(1)] Time- or bus-independent control: for any demand bus $v$, $(c^r_v, b^r_v, s^r_v) = (c_v, b_v, s_v)$ for all $1 \le r < R$ and some triple $(c_v, b_v, s_v)$;  or, for every demand bus $v$, 
$(c^r_v, b^r_v, s^r_v) = (c^r, b^r, s^r)$ for each $1 \le r < R$, for 
a certain triple $(c^r, b^r, s^r)$.
\item [(2)] Time-dependent, componentwise control.  This is a control such that 
for any given component $K$ of $G^r$, $(c^r_v, b^r_v, s^r_v)$ equals a fixed
triple $(c^r_K, b^r_K, s^r_K)$ for every $v \in K$.  
\item [(3)] Segmented control.  Let $(\Sigma_1, \Sigma_2, \ldots, \Sigma_H)$ be a
partition of the demand buses.  Then we insist that for any round $r$,
$(c^r_v, b^r_v, s^r_v)$ takes a common value for all demand buses in a given
set $\Sigma_i$.
\end{itemize}
An example of type (3) is that where the $\Sigma_i$ are quantiles of the demand distribution.  The resulting control is ``fair'' in that demands of 
similar magnitudes are reduced by similar fractions.  Controls of type (2),
in the case of the deterministic outage rule (F.1) can be explicitly 
described a-priori in polynomial space. This follows because in any fixed round
$r$, $G^r$ will have at most $n$ components and these depend solely on the
structure of the control on rounds up to $r-1$.  Thus in total at most $Rn$
distinct triples $(c^r_K, b^r_K, s^r_K)$ need to be specified.

Our primary focus are on algorithms and implementations for the most general version (time- and bus-dependent controls) of our approach, using the outage
rule (\ref{actual1})-(\ref{actual3}) with memory.  This is given in
Section \ref{firstorder}.  In Section \ref{scalingproblem} we will discuss a 
special case where the optimal control can be efficiently computed.

\section{First set of experiments}\label{firstset}
To motivate our overall approach, we first present experimental results using special, simple cases of the control.  The objective of these experiments is 
to expose, first, the fact that in cases of severe contingencies 
the application of control so as to immediately stop a cascade may be myopic
and could result in very suboptimal results.   In other words, it may pay 
off to allow some lines to become outaged.  However, this clashes with the
intuitive notion that postponing action to much later in the cascade results
in increasing uncertainty (because from the perspective of an agent at the
start of the cascade, the later stages of the cascade should be more 
uncertain).  We want to measure the impact of uncertainty 
in the timing of control, and to contrast it with the goal of avoiding
myopic, immediate action, as described above.

\subsubsection{The data}
In these experiments we used a snapshot of the U.S. Eastern Interconnect, with 
approximately $15000$ buses, $23000$ lines, $2000$ generators and $6000$
load buses.  The snapshot includes generator output levels, demands, 
and line parameters. 

Approximately $5000$ of the line flow limits were zero,
likely indicating a data error or missing data -- the minimum nonzero line 
flow value was $3 \times 10^{-2}$. When the flow limit 
of a line $j$ was equal to zero, we proceeded as follows, 
where $f^0_j$  is the initial power flow value on line $j$:
\begin{itemize}
\item If $|f^0_j| \ge 10^{-6}$, we reset the flow limit to  
$(1 + \gamma)|f^0_j|$, where $\gamma = 0.2$.
\item Otherwise, we reset the flow limit to $\beta = 10^{-4}$.
\end{itemize} 
A 
very small
number of lines had positive and large capacity, but nevertheless the 
$|f^0_j|$ values were close (or identical) to the line flow limits; in
which case we increased the line flow limit by $25\%$.

The reader may wonder about the impact of these numerical choices.  Based on
our limited testing with different choices of $\gamma$ and $\beta$, 
the impact is very minor in terms of $\beta$. The same holds for $\gamma$ 
unless we choose much values (on the order of $10^2$) and even 
then it is more a matter of degree than structure in the
cascades we will describe below.

Approximately $250$ of the lines had negative reactance values.  While there
are valid reasons for the use of negative reactances, in the experiments below
we assumed data error and replaced each negative value by its absolute value.

\subsubsection{Methodology}\label{firstsetmethod}

In 
all the experiments the same approach was employed: first, we interdicted
the grid according to a synthetic contingency, then we computed our affine
control, and finally we studied the behavior of this control.

To obtain contingencies we used the following methodology, which removes a
a set of $K$ random, high power flow lines from the grid, while preserving
connectivity.  Here $K$ is a given, small integer, and as before $m$ is
the number of lines.
\begin{itemize}
\item [(1)] A spanning tree $T$ is computed.
\item [(2)] Let $\hat f$ denote the power flow vector corresponding to the
given demands and generator outputs.  Renumber so that $|\hat f_1| \ge |\hat f_1| \ldots \ge |\hat f_m|$.
\item [(3)] Let $0 < \pi < 1$. Run steps (a) and (b), initialized with $S = \emptyset$, until stopping in (b): \\
\noindent For $j = 1, \ldots, m$, \\
\hspace*{.1in}{\bf (a)} If line $j \notin T$,  then \\ \hspace*{.25in} with probability $\pi$ reset $S \leftarrow S \cup {j}$.\\
\hspace*{.1in}{\bf (b)} If $|S| = K$, stop.
\item[(4)] The set $S$ of lines is removed from the network, producing network $G$ 
in our cascade template.
\end{itemize}
We used values of $K$ ranging from $1$ to $50$,
and for $\pi$ we used values ranging from $.1$ to $.5$.  

\subsubsection{The experiments}

We considered a case with $K = 2$ (two lines removed)
and $R = 20$ rounds. In the computation of 
the moving averages of line overloads (eq. (\ref{memory-1})) we used
$\alpha = 0.9$.  

First we consider the pure deterministic case of line outages, that is to
say we use line outage rule (\ref{stochout}) with $\epsilon_r = 0$
for all $r$.  If no control is applied, then at the end of round of round $20$ the 
{\em yield} (percentage of demand still being served) is $2.47 \%$.The
cascade is characterized by extremely high line overloads; see Table \ref{evolcomparo}  (we will discuss implications of this below).   
At the start of round 1, in fact, the maximum line overload is $40.96$, indicating that, likely, several lines with low flow limits are overloaded.

We computed the best control where
\begin{itemize}
\item [(i)] $c^r_v = b^r_v = 1$ for all $v$ and $r$.
\item [(ii)] $s^r_v = 0$ for all $v$ and $10 < r$. Thus, no control is applied
after round 10.
\item [(iii)] For each $1 \le r \le 10$, \underline{either} $s^r_v = 0.005$ for all
$v$, \underline{or} $s^r_v = 0$ for all $v$.
\end{itemize}
Thus we simply want to decide {\em when} to apply a control of a very simple
form.  Further, we are restricted to applying control in the first half of
the cascade; this is done as protection against uncertainty in the later
rounds of the cascade
The rationale for the numerical values in (iii) is that $1 + 0.005*(1 - 40) \approx 0.80$, that is to
say, the application of this control in round 1 will ``only'' shed $20\%$ of 
the demand.   

We want to stress that the experiments in this section do not amount to a rigorous attempt at optimizing control. In fact, the control obtained through
(i)-(iii) is only near-optimal.
Instead we are trying to provide an
example of the difference between an adequate control and the no-control option,
and the questions that arise from the comparison.  In particular, the amount $0.005$ was arrived at through a simple grid-search process.  

\begin{table}[h]
\centering
\caption{{\bf \emph{Cascade evolutions}}} 
\vskip 2 pt
\begin{tabular}{|c| r r r r| r r r r|}
\hline
\multicolumn{1}{|c|} {} & \multicolumn{4}{c|}{{\bf No control}} & \multicolumn{4}{c|}{{\bf c20}}\\
{\bf r} & \boldmath{$\kappa$} & {\bf O} & {\bf I} & {\bf Y}  & \boldmath{$\kappa$} & {\bf O} & {\bf I} & {\bf Y} \\\hline \hline
\noalign{\smallskip}
{\bf 1} & 40.96 & 86 & 1 & 100 & 40.96 & 86 & 1 & 100 \\  
{\bf 2} & 8.60 &187  & 8 & 99 &  8.60 & 165 & 8 & 96 \\  
{\bf 3} & 55.51  & 365 &20  & 98 & 61.74 & 303 & 17 & 96 \\  
{\bf 4} & 67.14  & 481 & 70  & 95 & 66.63 & 408 & 44 & 94 \\  
{\bf 5} & 94.61  & 692 & 149 & 93 & 131.08 & 492 & 94 & 93 \\  
{\bf 6} & 115.53  & 403 & 220  & 91 & 112.58 & 416 & 146 & 90 \\  
{\bf 7} & 66.12  & 336 & 333 & 89 & 99.62 & 326 & 191 & 78 \\  
{\bf 8} & 47.83  & 247 & 414 & 87 & 60.95 & 227 & 248 & 77 \\  
{\bf 9} & 7.16  & 160 & 457 & 85 & 32.50 & 72 & 279 & 76 \\  
{\bf 10} & 7.06  & 245 & 542 & 84 & 9.50 & 43 & 292 & 76 \\  
{\bf 11} & 37.55 & 195 & 606  & 83 & 45.28 & 35 & 303 & 76 \\  
{\bf 12} & 13.04 & 98 & 646 & 82 &  11.60 & 10 & 306 & 76 \\  
{\bf 13} & 22.61 & 128 & 688 & 82 & 3.88 & 6 & 310 & 75 \\  
{\bf 14} & 10.64   & 107 & 715 & 81 & 1.46 & 4 & 312 & 75 \\  
{\bf 15} & 5.03   & 64 & 721 & 81 & 1.34 & 1 & 312 & 75 \\  
{\bf 16} & 84.67 & 72 & 743  & 80 & 1.13 & 1 & 312 & 75 \\  
{\bf 17} & 32.15  & 52 &756 & 80 & 1.38 & 2 & 312 & 75 \\  
{\bf 18} & 6.50  & 43 & 763  & 80 & 1.26 & 1 & 312 & 75 \\  
{\bf 19} & 9.97  & 85 & 812 & 80 & 0.99 & 0 & 312 & 75 \\  
{\bf 20} & 32.34   & 39 & 812 & 2 & 0.99 & 0 & 312 & 75 \\  
\hline
\end{tabular}
\label{evolcomparo}
\end{table}

In any case, the optimal control that satisfies conditions (i)-(iii) 
(and which we shall refer to as {\bf c20} for future reference) attains
a termination yield of $75.2\%$, picks
rounds $2$ and $7$ to apply control.  Note that 
since the maximum line overload is high
in round $1$, c20 allows some lines to become outaged in round $1$.

This point is further elaborated in Table 1, where
 ``r'' indicates round and for each round, ``$\kappa$'' 
indicates maximum line overload at the start of the round, ``O'' is the number of lines outaged during the round, ``I''
is the number of islands at the end of the round and ``Y'' is 
the (rounded) percentage
of demand
being delivered end of the round.  {\bf We stress that we count {\em all} islands,
even those that consist of a single bus with no demand, and when computing the 
maximum line overload we consider {\em all} lines, no matter how minor.}\\

\noindent {\bf Discussion.} We see that initially both cascades have extremely high line overloads, many line outages and large amounts of islanding.  However, 
under c20 after line 11
the overages are significantly smaller, and rapidly
decreasing, and after round 8 the number of new outages and islands is also much
smaller (and decreasing); both in spite of the fact that control is 
last applied in round 7.  Thus, effectively, the cascade has been ``stabilized''
under c20, long before the end of the time horizon.

The reader might wonder about the rapid decrease of yield from $80\%$ to
$2\%$ in the no-control case.  This is due to the termination feature in
our cascades that requires all line overloads to be eliminated by the end
of the last round;  since the no-control cascade has very high
maximum overload ($32.34$),
at the start of round $20$, the termination rule forces a drastic reduction in yield.

Nevertheless, in the no-control case,
the combination of comparatively high yield (up to
round 4), high number of line
outages, large line overloads and large amount of islanding suggest the possibility that many of the 
outages involve unimportant lines, and likewise with many of the islands (though of course a $22\%$ yield loss should indicate a severe contingency).  One wonders if somehow the no-control option {\em might} be attractive
if {\em enough} time (i.e., rounds) were available. 

\begin{table}[h]
\centering
\caption{{\bf \emph{Further evolution of no-control cascade from Table \ref{evolcomparo}}}} 
\vskip 2pt
\begin{tabular}{|c | r r r r r r r r|}
\hline
{\bf r} & {\bf 25}  &  {\bf 28} & {\bf 29} & {\bf 30} & {\bf 31} & {\bf 32} & {\bf 33} & {\bf 34}   \\ 
\hline
\noalign{\smallskip}
{\bf O}  & 21.63  & 2.00 & 5.70 &2.50  & 2.38 & 1.35 & 1.07 & 0.99  \\
{\bf Y} & 79 & 78 & 78 &78  & 78 & 78 & 78 & 78  \\
\hline
\end{tabular}
\label{table2}
\end{table}

To investigate these possibilities, we extended the no-control cascade. Table
\ref{table2} shows the results for selected rounds.  We see that the no-control
approach finally yields stability by round 34, attaining yield $78\%$.  
This is slightly better (but very close) to what c20 obtained in 20
rounds (and, furthermore, control action under 
c20 was restricted to rounds 1-10).  
Nevertheless, the no-control approach experiences significant line overloads 
as late as round 32.

By maintaining high overloads into very late rounds, the no-control strategy
becomes more exposed to the unavoidable {\em uncertainty} that should be
taken into account when modeling cascades, and which we have up to now 
ignored.  We model noise by means of fault outage rule (\ref{stochout}).  In the following set of tests we assume that 
\begin{eqnarray}
\epsilon_r & = & 0.01 \, + 0.05* \lfloor r/10 \rfloor. \label{noise1}
\end{eqnarray}
Possibly, noise should be increasing at a faster rate than the 
above formula stipulates (perhaps exponentially).  However, the control
considered in Table \ref{evolcomparo} as well
as the no-control approach are both exposed to significant amounts of 
noise after round 10; more so in the no-control case.  
We would thus expect that under
rule (\ref{noise1}) the no-control approach will perform much more poorly.

To test these hypothesis, we ran 1000 simulations of cascades under rule (\ref{noise1}) for the no-control case and for control using c20.  The results are
summarized as follows: using c20, the average yield is $42.90$ and the standard
deviation of yield is $27.47$, whereas using no control the average yield 
is $7.96$ and the standard deviation is $9.33$.  In other words, c20 proves
much more robust than the no-control strategy, which is not surprising given
the structure of rule (\ref{noise1}).  A question that arises as a result is
whether c20 is in some sense optimally robust.

One way to investigate this question is to investigate controls that are
less exposed to uncertainty by restricting them to a shorter timeline, i.e.
by enforcing termination before round 20. For $T = 10, \, 15, \, 25$, we compute an optimal control required
to terminate by round $T$, and otherwise subject to rules 
(i)-(iii), that is 
$c^r_v = b^r_v = 1$ for all $v$ and $r$, $s^r_v = 0$ for all $v$ and $10 < r$,
and for each $1 \le r \le 10$, either $s^r_v = 0.005$ for all
$v$, or $s^r_v = 0$ for all $v$.  We name these controls c10, c15 and c25, respectively.

Table \ref{table3} presents the comparisons between all the options we 
have considered. In this table, ``DetY'' is the yield in the deterministic
case ($\epsilon_r = 0$ for all $r$), ``MaxY'' and ``MinY'' are the maximum 
and minimum yields in all
the simulations (resp.), ``AveY'' is the average yield and ``StddY'' is the standard
deviation of yield.

\begin{table}[h]
\centering
\caption{{\bf \emph{Robustness comparison - 1000 runs using stochastic outage rule (\ref{stochout}) with noise as in (\ref{noise1}) }}} 
\vskip 2pt
\begin{tabular}{|c | r r r r r|}
\hline
{\bf Option} & {\bf DetY} & {\bf MaxY}  &  {\bf MinY} & {\bf AveY} & {\bf StddY} \\ 
\hline
\noalign{\smallskip}
{\bf c10} & 37.49 & 39.05 & 0.00 & 11.81 & 11.84 \\
{\bf c15} & 72.44 & 71.85 & 0.00 & 33.94 & 22.51 \\
{\bf c20} & 75.19 & 76.30 & 1.17 & 41.90 & 27.47 \\
{\bf c25} & 77.23 & 42.34 & 1.38 & 11.99 & 10.97 \\
{\bf no control} & 77.75 & 36.04 & 0.00 & 7.96 & 9.33 \\
\hline
\end{tabular}
\vskip -5pt
\label{table3}
\end{table}

Control c20 emerges as superior over c15 and c10.  This can be
explained as follows. Even though c15 and c10
are significantly less exposed to risk than c20, 
they are also restricted to
operating, and terminating,
during a stage of the cascade characterized by extremely high line overloads.
Control c20, by being able to operate over 20 rounds, has ``more time'' 
while also avoiding the large uncertainty rounds 20 and higher.
For this reason, c20 is also
superior to c25 (their averages are separated by more than one standard 
deviation).
One common feature
that emerges in controls c10, c15, c20 and c25 (not shown in the table) is that no control is taken in 
round 1, and control is taken in round 2 (and in the cases of c10, c15 and c20, rounds 5 or 7).

We stress that (\ref{noise1}) is {\em one} categorization of noise.  Using a different formula the outcome could be different, say c15 could prove best.  However, the outlook we are taking here is that by computing {\em a} robust control
with respect to {\em some} rule such as (\ref{noise1}) we obtain a control
that remains robust (though possibly not optimally so) even if the model for
uncertainty were to be somewhat changed.  And, in any case, computing a control
which is is somewhat robust should be better than completely ignoring uncertainty.  

To explore these issues, we study the following model
\begin{eqnarray}
\epsilon_r & = & 0.01 \, + 0.005* r, \label{noise2}
\end{eqnarray}
which can be considered a smoothed version of (\ref{noise1}).  Under this
model both c15 and c20 are exposed to more noise than c10, and more noise than under rule (\ref{noise1}).  Consider Table \ref{table4}.
\begin{table}[h]
\centering
\caption{{\bf \emph{Robustness comparison - 1000 runs using stochastic outage rule (\ref{stochout}) with noise as in (\ref{noise2}) }}} 
\vskip -5pt
\begin{tabular}{|c | r r r r r|}
\hline
{\bf Option} & {\bf DetY} & {\bf MaxY}  &  {\bf MinY} & {\bf AveY} & {\bf StddY} \\ 
\hline
\noalign{\smallskip}
{\bf c10} & 37.49 & 38.93 & 0.00 & 7.54 & 9.55 \\
{\bf c15} & 72.44 & 63.94 & 3.41 & 28.02 & 17.94 \\
{\bf c20} & 75.19 & 73.04 & 0.00 & 32.24 & 21.30 \\
{\bf c25} & 77.23 & 54.62 & 0.25 & 16.84 & 12.66 \\
{\bf no control} & 77.75 & 18.86 & 0.00 & 5.11 & 5.28 \\
\hline
\end{tabular}
\label{table4}
\end{table}
We see that c20 still appears superior to the other controls, though c15
is almost as good.\\

The above experiments do not amount to a full optimal
robust control computation. In Section \ref{gridstoch} we will return to these experiments from a stochastic optimization 
perspective.

\section{Optimization methods}\label{methodology}
Given a control vector $(c,b,s)$, denote by $\tilde \Theta^R(c,b,s)$ the 
final demand at termination of the $R$-round cascade controlled by $(c,b,s)$.
Our goal is to maximize $\tilde \Theta^R(c,b,s)$ over all controls.  This is
a nonconcave, in fact very combinatorial, maximization problem \cite{boco}, \cite{luen}; it is very large (e.g. if $R = 10$ the $(c,b,s)$ vector has more than
$180000$ variables in the case of the Eastern Interconnect).
 It is also important to incorporate stochastics.  

In principle, the deterministic
case of our problem could be tackled using mixed-integer programming techniques,
and the stochastic version, using stochastic programming
\cite{lind}.  Of course, one could choose a different formulation of the
cascade control problem than the one we chose (using a different kind of
control, for example).  But any formulation will have to deal with 
the combination of combinatorics in the network dynamics, multistage
behavior,
stochastics and very large size. In our opinion, this combination  
places the problem outside the capabilities of current optimization 
methodology, even in the deterministic case.  We remind the reader that 
we envision our control as being computed in real time and 
we might only have one hour, or less, to do so.  

Another point to stress is that nonconcavity in a maximization problem leads to
non-monotone behavior: in our case, just because a small change in control
leads to an improvement does not imply that a larger change will result in
greater improvement.

\subsection{First-order methods for the general case}\label{firstorder}
Here we describe a procedure to compute a control given by triples $(c^r_v, b^r_v, s^r_v)$, for each demand bus $v$ and each round $r$, and using a control law as in 
Step 2 of algorithm (\ref{affinecontcasc}). 
We will assume a semi-random outage rule as in (\ref{actual1})-(\ref{actual3}), 
with memory, as in (\ref{memory-1}).  As previously, the goal is to compute a control
that will maximize the expectation of the amount of demand being served by the end of round $R$, which as before
we denote by $\tilde \Theta^R(c,b,s)$.  We stress that the control parameters we use 
are state-independent; this is a design feature.

Toward this goal we will use an algorithm based on the following template:
\begin{center}
  \fbox{
    \begin{minipage}{0.9\linewidth}
      \hspace*{.9in} \begin{PRO}{First-order algorithm}\label{stochgrad1} \end{PRO}
      {\bf Input}: a control vector $(c,b,s)$.\\
      \For $k = 1, 2, \ldots$ \Do \\
        \hspace*{.5in}{\bf 1.} Estimate $g = \nabla \tilde \Theta^R(c,b,s)$.\\
        \hspace*{.5in}{\bf 2.} Choose ``step-size'' $\mu \ge 0$ and update control
        to $(c,b,s) \, + \, \mu (g_c, g_b, g_s)$.\\
        \hspace*{.5in}{\bf 3.} If $\mu$ is small enough, stop.
    \end{minipage}
  }
\end{center}

This is a common first-order (steepest-ascent) method.  
In the deterministic case, Step 1 should be interpreted as a an approximate
rule since $\tilde \Theta^R$ is not differentiable (our stochastic outage rule
\ref{stochout} does smooth  out the expectation).  The vices of procedure
\ref{stochgrad1} are well known: even if $\tilde \Theta^R$ were smooth, its 
nonconcavity implies that the steepest-ascent method may not converge to
a global optimum.  And even if $\tilde \Theta^R$ were smooth and concave, 
steepest ascent may zigzag or stall.  See \cite{luen}.  

In summary, 
Procedure \ref{stochgrad1} should be viewed as a {\em local search} method
with which to explore the neighborhood of a solution.  Finally, in our 
setting the procedure could prove expensive, since each evaluation of $\tilde \Theta^R$
(including in the estimation of $\nabla \tilde \Theta^R$ through finite 
differences) requires a cascade simulation, each round of which requires two
power flow computations in our setup.

On the positive side, however, the procedure is flexible enough to handle
(at increased computational cost) important features, such as more realistic
AC power flow models, or more complete renditions of low-level controls in
the operation of a power grid.  Essentially, Procedure \ref{stochgrad1}
is an example of simulation-based optimization, i.e. it
only needs to have a ``black-box'' that computes the function $\tilde \Theta^R$.

An active research field that considers optimization under such assumptions
is that of {\em derivative-free optimization} (see \cite{derivfree}) 
and related methods that incorporate second-order information \cite{ipopt}.  
In our estimation, these methodologies may not scale well to problems of 
the size we consider.   In forthcoming
work we will experiment on adaptations of these methodologies to our 
problem.

When we consider a model that includes stochastics, the first-order method 
can be viewed as a {\em stochastic gradients} algorithm (see  \cite{robbinsmonro}, \cite{kushnerclark} -- an alternative methodology is
provided by bundle methods).  In the stochastic gradients approach, a fixed {\em sample path} of the appropriate random variables is chosen in advance of
each gradient and step-length computation.  In 
Section \ref{stochasticoptimization} we will further discuss this approach. 

Whether we use the stochastic setting or not, we cannot completely rely on
Procedure \ref{stochgrad1} as the sole optimization engine -- to repeat the
above,  the resulting algorithm
would both be too slow and likely to get trapped in local 
maxima.  To help avoid these difficulties we rely on several heuristics
described later.  In the next section we describe a special case of the
optimal control problem that can be efficiently solved.

\subsection{The optimal scaling problem}\label{scalingproblem}
In this section we describe an algorithm that computes an optimal time-dependent componentwise control under outage rule (F.1), without memory. Either version of rule (F.1) can be used; for simplicity of language we will use the first.
For brevity, we will refer to this as the {\em simple scaling setting}.
Our algorithm computes an optimal control
in time $O(m^{R-1}/(R-1)!)$ where as before $m$ is the number of lines.  

\begin{RE} \label{triplestosingle} Consider an optimal control. 
Let $1 \le r < R$ and let $K$ be a component of $G^r$ under the optimal 
control.  Then (at round $r$) we
will scale all demands in $K$ by a common multiplier $0 \le \lambda^r_K \le 1$ (defined
as in Procedure \ref{affinecontcasc}.  Clearly, the control can
be equivalently defined by the values $\lambda^r_K \, (\le 1)$ rather than the triples
$(c^r_K, b^r_K, s^r_K)$, and we will use this convention below.\end{RE}

\begin{NO} Let $G$ be a graph, and let $\mu$ be a supply-demand vector on
$G$.  We denote by $\hat f(G, \mu)$ the unique, feasible flow vector on
$G$ when $\mu$ is the supply-demand vector (see Remark \ref{uniqueness}).
\end{NO}

\noindent In what follows we assume that we have a given
supply-demand vector $\beta$.  Let $R$ be the number of
rounds for the cascade.  Our problem is to compute a control that maximizes 
the total demand satisfied after $R$ rounds, assuming that at the start of
round 1, $\beta$ is the supply-demand vector.  We will solve this as a 
special case of a family of problems:  
\begin{DE} For $t \ge 0$ real, denote by
$\Theta^{(R)}_G(t | \beta)$ 
the final total demand resulting from applying an optimal 
control in an $R$-round cascade on graph $G$, where the initial supply-demand vector
is $t \beta$.  
\end{DE}
We will show that $\Theta_G^{(R)}(t | \beta)$ is a nondecreasing
piecewise linear function of $t$ with at most $Rn$ pieces.

\begin{RE} \label{scalingflows} Let $\alpha$ be a supply-demand vector
on graph $G$.  Let $t \ge 0$. Then 
$\hat f(G, t \alpha) = t \hat f(G, \alpha)$.
\end{RE}

\begin{LE}\label{lemma1}
Let $\mu$ be a supply-demand vector.  Suppose $G$ is connected. 
Then 
$\Theta_G^{(1)}(t | \mu)$ is a nondecreasing piecewise-linear function of $t$ with two pieces.
\end{LE}
{\em Proof.} Note that since $R = 1$, only Steps 1 and 2 in algorithm
(\ref{affinecontcasc}) will be executed.  Further, 
writing $\hat f = \hat f(G^1, \alpha)$, 
when running (\ref{affinecontcasc}) starting
with the initial supply-demand vector $t \mu$, we will have
$f^1 = t \hat f$ in Step 1, and  writing $\psi = \max_{j} |\hat f_{j} |/u_{j}$, we have that $\max_{j} |f^1_{j} |/u_{j} \, = \, t \psi$.
Denoting  by $\tilde D$ the sum of demands implied by $\mu$ we have
as per our cascade termination criterion that the final total demand at the end of
$R = 1$ rounds will equal
\begin{eqnarray}
t \tilde D, && \mbox{if} \ \  t \, \le \, 1/\psi, \ \ \ \mbox{and} \\
\frac{t}{t \psi} ~ \tilde D \ = \frac{1}{\psi}~ \tilde D, && \mbox{otherwise}.\ \ \hspace*{.2in}\qed 
\end{eqnarray}
Now we turn to the general case with $R > 1$.  We assume, without loss of
generality, that $G^1$ is connected.  Let $\hat f = \hat f(G^1, \beta)$.
\begin{DE} A {\em critical point} is a real $\gamma > 0$, such that for
some line $j$, $\gamma \hat f_{j} = u_{j}$.
\end{DE}
Recall that we assume $u_j > 0$ for all $j$; thus let $0 < \gamma_1 < \gamma_2 < \ldots < \gamma_p$ be the set of all distinct
critical points.  Here $0 \le p \le m$.  Write $\gamma_0 = 0$ and 
$\gamma_{p+1} = +\infty$. 

\begin{DE}
For $1 \le i \le p$ let $F^i = \{ j \, \in \, \cA \, : \, \gamma_h | \hat f_{j} | = u_{j} \}$. 
\end{DE}
\noindent Now assume 
that the initial supply-demand vector is $t \beta$ with $t > 0$ and
let $0 < \lambda^1 \le 1$ be the optimal multiplier used to scale demands
in round 1 (see Remark \ref{triplestosingle}).  Write 
\begin{eqnarray}
 && q = \argmax\{ h \, : \, \gamma_h < t \}. \label{defq}
\end{eqnarray}
Thus, $t \le \gamma_{q+1}$, and so $\lambda^1 t \le \gamma_{q+1}$. We stress that 
these relationships remain valid in the boundary cases $q = 0$ and $q = p$.

\begin{NO} Let the index $i$ be such that $\lambda^1 t \in (\gamma_{i-1}, \gamma_{i}]$. \end{NO}

\noindent Note
that in Step 3 of round 1 
of algorithm (\ref{affinecontcasc}) we will scale 
all demands by $\lambda^1$, and since we assume $G^1$ is connected, 
in Step 4 we will also scale all supplies by 
$\lambda^1$. Thus, for any
$h \le i-1$, and any line $j \in F^h$, we have that after Step 4 
the absolute value of the flow on $j$ 
is 
\begin{eqnarray}
\lambda^1 \,t \, | \hat f_{j} | > \gamma_h \, | \hat f_{j} | = u_{j},
\end{eqnarray}
and consequently $j$ becomes outaged in round 1.  On the other hand,
for any line $j \notin \cup_{h \le i-1} F^h$, the absolute value of the
flow on $j$ immediately after Step 4 is
\begin{eqnarray}
\lambda^1 \,t \, | \hat f_{j} | \le \gamma_{i} \, | \hat f_{j} | \le u_{j},
\end{eqnarray}
and so $j$ does not become outaged in round 1.  In summary, the
set of outaged lines is $\cup_{h \le i-1} F^h$;  in other
words, we obtain the same network $G^2 = G^1 \, - \, \cup_{h = 1}^{i-1} F^h$ for every
$t$ with $\lambda^1 t  \in (\gamma_{i-1}, \gamma_{i}]$.  

\begin{NO} For an index $j$, write $\cK(j)$ = set of components of $G^1 \, - \, \cup_{h = 1}^{j} F^h$. \end{NO}

Let $H \in \cK(i-1)$.  Then, prior to Step 6 of round 1, the supply-demand vector for $H$
is precisely the restriction of $\lambda^1 t \beta$ to the buses of $H$, and when
we adjust supplies and demands in Step 6, we will proceed as follows 
(where we use notation as in Section \ref{notation}):
\begin{itemize}
\item if $ \sum_{s \in \cD \cap H} (-\lambda^1 t \beta_s) \, \ge \, \sum_{s \in \cG \cap H} (\lambda^1  t \beta_s)$ then for each demand bus $s \in \cD \cap H$ we will reset
its demand to
$$ -r \lambda^1  t \beta_s, \ \ \mbox{where} \ \ r = \frac{\sum_{s \in \cG \cap H} (\lambda^1 t \beta_s)}{\sum_{s \in \cD \cap H} (-\lambda^1  t \beta_s)} = -\frac{\sum_{s \in \cG \cap H} (\beta_s)}{\sum_{s \in \cD \cap H} (\beta_s)},$$
and we will leave all supplies in $H$ unchanged. 
\item likewise, if
$ \sum_{s \in \cD \cap H} (-t \lambda^1 \beta_s) \, < \, \sum_{s \in \cG \cap H} (\lambda^1  t \beta_s)$ then the supply at each bus $s \in \cG \cap H$ will be reset
to
$$ r \lambda^1 t \beta_s, \ \ \mbox{where} \ \ r = -\frac{\sum_{s \in \cD \cap H} (\beta_s)}{\sum_{s \in \cG \cap H}} (\beta_s),$$
but we will leave all demands in $H$ unchanged.
\end{itemize}
Note that in either case, in round 2 component $H$ will have a supply-demand vector of
the form $\lambda^1 t \beta^H$, where $\beta^H$ is a supply-demand vector.  
Thus an optimal control on $H$, on rounds $2, \ldots, R$, 
will yield a final total demand 
\begin{eqnarray}
&& \Theta_H^{(R-1)}(\lambda_1 t | \hat \beta^H), \label{Hfinal}
\end{eqnarray}
which, inductively, is a nondecreasing function of $\lambda_1 t$, and therefore
is largest when 
\begin{eqnarray}
&& \lambda^1 \ = \ \min\left\{1, \frac{\gamma_{i}}{t} \right\}.
\end{eqnarray}

\vspace{.1in}
\noindent {\bf Case 1.} Suppose $i \le q$. As noted above, by definition
(\ref{defq})
of $q$  we have that $\gamma_{i} \le \gamma_{q} < t$.  Thus, the expression in (\ref{Hfinal})
is maximized when $\lambda_1 = \frac{\gamma_{i}}{t}$,
and we obtain final ($R$-round)
demand equal to
\begin{eqnarray}
&& \sum_{H \in \cK(i-1)} \Theta_H^{(R-1)}(\gamma_{i} | \hat \beta^H), \label{fixedcase1}
\end{eqnarray}
which is independent of $t$. 

\vspace{.1in}
\noindent {\bf Case 2.} Here $q < i$, and so $i = q+1$ by definition of $q$
and $\lambda^1 \le 1$. Thus (\ref{Hfinal}) is maximized by setting $\lambda^1 = 1$. The final demand equals
\begin{eqnarray}
&& \sum_{H \in \cK(q)} \Theta_H^{(R-1)}(t | \hat \beta^H). \label{case2}
\end{eqnarray}

\noindent In summary, we have:
\begin{eqnarray}
 \Theta^{(R)}_G(t|\hat \beta) \ = \ \max \left\{ \, \max_{1 \le i \le q} \left\{\sum_{H \in \cK(i-1)} \Theta_H^{(R-1)}(\gamma_{i} | \hat \beta^H)\right\}  \ , \ \sum_{H \in \cK(q)} \Theta_H^{(R-1)}(t | \hat \beta^H) \, \right\}. \label{consolidated}
\end{eqnarray}
\begin{THM}\label{th2}
(i) $ \Theta^{(R)}_G(t|\hat \beta)$ is nondecreasing, piecewise-linear, 
with at most 
$$ \frac{m^{R-1}} {(R-1)!} ~ + ~ O\left(m^{\max\{1, R-2\}}\right)$$ 
breakpoints.\\
\noindent (ii) The optimal choice for $\lambda^1$ is $\lambda^1 = 1$ or 
$\lambda^1 = \gamma_k/ t$ for some $k$.
\end{THM}
\noindent {\em Proof.} (i) By induction on $R$, starting from Lemma \ref{lemma1}.  For the general step, consider the above discussion which assumes that
$\lambda^1 t \in (\gamma_{i-1}, \gamma_{i}]$.  Then if Case 1 above holds,
we have that $ \Theta^{(R)}_G(t|\hat \beta)$ is constant.  And if 
Case 2 holds, then equation (\ref{consolidated}) applies.  The form of 
(\ref{consolidated}) guarantees that, inductively, 
$ \Theta^{(R)}_G(t|\hat \beta)$ is nondecreasing piecewise-linear. 

To analyze the number of breakpoints in $\Theta^{(R)}_G$, assume
first that $R = 2$.  Consider the
effect of removing, {\em  one at a time}, the lines of  $\cup_{h = 1}^{p} F^h$.
Prior to its removal, each line $j$ has both ends in the same component $K$; 
the
removal either creates two new components (if $j$ is a bridge of
$K$) or creates a new component (which differs from $K$ in that line $j$ is
not included).  Thus the removal process can be represented as a binary tree
whose leaves correspond to the components of $G^1 \, - \, \cup_{h = 1}^{p} F^h$,
i.e. the members of $\cK(p)$.  Since these are disjoint there are at most
$n$ of them; since in a binary tree the number of degree three vertices is
at most the number of leaves we conclude that 
$$ | \cup_{h = 1}^p \cK(h) | \le m  + n \le 2 m + 2.$$
Furthermore, let $H$ be a component in
$ \cup_{h = 1}^p \cK(h)$. Define $h = \min \{ j \, : \, H \in \cK(j) \}$
and $h' = \max \{ j \, : \, H \in \cK(j) \}$.  By Lemma \ref{lemma1}, it
follows that 
$\Theta^{(1)}_H$ will contribute at most one breakpoint
to $\Theta^{(R)}_G$, and 
that this breakpoint will occur for some $t$ with $\lambda^1 t \in [\, \gamma_h  \, , \, \gamma_{h'} \,)$. The maximum in 
(\ref{consolidated}) shows that for each $q$, one additional new breakpoint
is created.  Thus, in total, $\Theta^{(R)}_G$ has at most $O(m)$ 
breakpoints and the result
is verified for $R = 2$.\\

\vspace{.1in}

In what follows we assume that $R \ge 3$.  Suppose
$q = 0$ and thus $i = 1$.  Since $\lambda_1 t < \gamma_1$, it follows that
no lines are outaged in round 1, i.e. $G^2 = G^1 = G$, and in subsequent rounds
no line will be overloaded.  Thus, in this case, $\Theta^{(R)}_G(t|\hat \beta) = t \tilde D$ and there are no breakpoints. For $q > 0$ we proceed using (\ref{consolidated}). For each $H \in \cK(q)$, 
inductively,
$\Theta^{(R-1)}_H$ has at most 
$$ \frac{m_H^{R-2}}{(R-2)!} \ + \ c \, m_H^{\max\{1, R-3\}} $$ 
breakpoints, where 
$m_H$ denotes the number of lines in $H$ and $c \ge 0$ is a constant.
So (\ref{consolidated}) implies that subject to $i = q+1$, the
number of breakpoints in $\Theta^{R}_G$ 
is at most
\begin{eqnarray}
&& ~ 1 ~ + ~ \sum_{H \in \cK(q)}\left[ \frac{m_H^{R-2}}{(R-2)!} \, + \, c \, m_H^{\max\{1, R-3\}}\right] \nonumber \\
&& \le \ 1 ~ + ~ \frac{\left( \, m \, - \, | \cup_{h = 1}^{q} F^h| \right)^{R-2}}{(R-2)!} ~ + ~ c \left( \, m \, - \, | \cup_{h = 1}^{q} F^h| \right)^{\max\{1,R-3\}} \nonumber \\
&& \le 1 ~ + ~ \frac{\left( \, m \, - \, q \right)^{R-2}}{(R-2)!} ~ + ~ c \left( \, m \, - \, q \, \right)^{\max\{1,R-3\}}
\label{sumconsolidated}
\end{eqnarray}
since  $\left( \cup_{h = 1}^{q} F^h \right) \cap H = \emptyset$ for each 
$H\in \cK(q)$.  Summing this expression over all $1 \le q \le p$, we
obtain that the total number of breakpoints is at most
\begin{eqnarray}
&&  p ~ + ~ \sum_{q = 1}^p\left[ \frac{\left( \, m \, - \, q \right)^{R-2}}{(R-2)!} ~ + ~ c \left( \, m \, - \, q \right)^{\max\{1,R-3\}} \right] \nonumber \\
&& \le m ~ + ~ \sum_{q = 1}^m\left[ \frac{\left( \, m \, - \, q \right)^{R-2}}{(R-2)!} ~ + ~ c \left( \, m \, - \, q \right)^{\max\{1,R-3\}} \right] \nonumber \\
&& \le \frac{\left( \, m \, - \, 1 \right)^{R-1}}{(R-1)!} ~ + ~ O((m-1)^{R-2}) + ~ m ~ + ~ c \sum_{q = 1}^m \left( \, m \, - \, q \right)^{\max\{1,R-3\}}. \label{hew}
\end{eqnarray}
For $R = 3$ the last three terms in (\ref{hew}) are $O(m)$ and we are done as desired.  For 
$R > 3$, the last term in (\ref{hew}) equals
\begin{eqnarray}
&&  c \frac{\left( \, m \, - \, 1 \right)^{R-2}}{R - 2} ~ + ~ O(m^{R-3}),
\end{eqnarray}
and again we conclude as desired for $c $ large enough.\\

\noindent (ii) This follows from the discussion leading to eq. (\ref{consolidated}). \QED

\vspace{.1in}
Part (ii) of Theorem \ref{th2} illustrates a weakness of the simple scaling
approach -- when applying an optimal control, at least one line becomes
fully loaded at each round.  Such a strategy is likely non-robust.  We 
plan to address this issue in upcoming work;  using the stochastic
outage (F.2) and computing an appropriate optimal control.

Despite the apparent shortcomings of the method, and of the simplicity of
the proposed control, the ability to compute a global optimum in polynomial
time (for fixed $R$) is a significant asset, especially as a starting point
for the simulation-based methods for the general problem that are 
proposed below. 
In forthcoming work we will 
implement an appropriate version of the above algorithm; an relevant question
is whether the worst-case bound in Theorem \ref{th2} is attained using
realistic data.

\section{The algorithm}\label{thealgorithm}

Our algorithm implements Procedure \ref{stochgrad1} to implement 
an affine control as in Template (\ref{affinecontcasc}),
repeated here for convenience.
The control specifies, for each round $r$ of the cascade and each demand
bus $v$, a triple $(c^r_v, b^r_v, s^r_v)$. At round $r < R$ of the cascade, each
demand bus $v$ observes the maximum line load $\kappa^r_v$ in the 
component that $v$ currently belongs to. Then, where $d^r_v$ denotes the
current value of the demand at $v$,
\begin{eqnarray}
&& \mbox{if} \ \kappa^r_v > c^r_v, \ \mbox{demand at $v$ is reset to} \  \min\{1, \, [b^r_v + s^r_v (c^r_v - \kappa^r_v)]^+ \}~d^r_v. \label{reminder}
\end{eqnarray}
The ``normal'' case of such a control is that where $b^r_v = c^r_v = 1$, and
$s^r_v \ge 0$, which decreases demands in proportion to the maximum overload.
However, cases with $c^r_v \neq 1$ (delayed or proactive control) can
prove optimal.  Setting $b^r_v < 1$ can result in nonsmooth (fixed-penalty)
controls.  Finally, setting $s^r_v < 0$, though counterintuitive, can prove
optimal when non-monotone behavior occurs. \\

\noindent In order to initiate the gradient search method, we rely on 
grid-search, a standard enumerative idea:

\vspace{.1in}
\noindent {\bf Grid search.}  Here we fix $c^r_v = b^r_v = 1$ for all $r$ and 
$v$, and $s^r_v = 0$ for all $v$ and all $2 < r < R$.  Thus
the only remaining parameters are $s^r_v$ for all $v$ and $r =1, 2$.
We restrict the search to two values $\bar s^1$ and $\bar s^2$, 
and insist that for all $v$ and $1 \le r \le 2$ we 
have $s^r_v = \bar s^r$.   In our current implementation,
this two-dimensional search,
in turn, is carried out one parameter at a time, as follows.  Let $\tilde \kappa^1$
be the maximum line overload observed in the no-control cascade, during round
1.  Assuming $\tilde \kappa^1 > 1$, then we enumerate all choices for
$\bar s^1$ of the form $\bar s^1 = (.1 + .008\, i)/(\tilde \kappa^1 - 1)$ for
$i = 0, 1, \ldots, 100$.  In other words, this enumerates all controls where
in round 1 we scale demands by a factor of $.9 - .008 \,i$ for $i = 0, 1, \ldots, 100$.  Let $\bar s^1_1 < \bar s^1_2$ be the two enumerated choices which 
produce the highest and second highest $\tilde \Theta^R$ value.  Then we
repeat the search in the interval $[\bar s^1_1 , \bar s^1_2]$ by enumerating 
$\bar s^1 = \bar s^1_1 + i \, (\bar s^1_2 - \bar s^1_1)/100$ for $i = 0, 1, \ldots, 100$.  
The value that produces that highest
$\Theta^R$ value is our final choice for $\bar s^1$.  We fix this value and
now carry out the same type of search for $\bar s^2$.\\

\vspace{.1in}

We will see below that grid search can produce very good control vectors, but
which in general can be improved, sometimes significantly, by widening the
search.  One can use the control computed by grid search to start the general
gradient search; however in high-dimensional cases even general gradient search
itself can be quite slow as each gradient estimation step could prove very
slow.  This will not be the case if enough parallel computing resources
are available; however and we have found an additional step to be useful:
\vspace{.1in}

\noindent {\bf Segmented search}.  As introduced in Section \ref{adaptive1}, 
consider 
a fixed partition
$(\Sigma_1, \Sigma_2, \ldots, \Sigma_H)$ of the demand buses. We search,
using the first-order method,
for triples of the form $(\hat c^r_i, \hat b^r_i, \hat s^r_i)$ for 
each $1 \le r < R$ and $1 \le i \le H$, so as to obtain the control where
for each $1 \le r < R$,  and each demand bus $v$, 
$(c^r_v, b^r_v, s^r_v) = (\hat c^r_i, \hat b^r_i, \hat s^r_i)$ if $v \in \Sigma_i$.  In our implementation, the $\Sigma_i$ are {\em demand quantiles}. That
is to say, if $L$ is the number of demand buses, then $\Sigma_1$ contains
the $\lfloor H/L \rfloor$ buses with largest demand, $\Sigma_2$ 
contains the next $\lfloor H/L \rfloor$ buses with largest demand, and so on.
The advantage of this approach is that it considerably reduces the dimensionality of the problem, even if $H$ is chosen relatively large, such as
$H = 100$.  In fact, the (segmented) first-order method
runs quite fast, and the approach in \cite{ipopt} might also be practicable. 
Further, a segmented control is arguably 'fair' in that 
it specifies, to some degree, 
that similar buses are bound by similar control laws, though we stress that
when applying the control (\ref{reminder}) the actual demand reduction can
be very different for two buses in the same segment but in different components.

In the first set of experiments we have conducted, we have chosen $H = 50$ and
we fix $\hat b^r_i = 1$ for $1 \le r < R$ and $1 \le i \le H$.  Thus, altogether, we have $2 \, H \, (R - 1) = 100 \, (R - 1)$ variables, still large but
much more manageable than full gradient search.  In the experiments below,
we {\em forgo} full gradient search; however it would be straightforward to
follow up segmented search with full gradient search.

\subsection{Implementation details}\label{implementation}

The parallel implementation of our algorithm relies 
on the familiar master-worker paradigm.  Each worker performs computations of
the function $\tilde \Theta^R(c,b,s)$ for a given control $(c,b,s)$ whereas
the master carries out the gradient search algorithm.  In the experiments
we report on here, we use the linear power flow model; linear programs are
solved using Cplex 12.0 \cite{cplex} and Gurobi 3.0 \cite{gurobi}.  These solvers were used with all presolve options turned-off (this increased robustness).
Further, the flow component in the solution to 
the linearized power flow system (\ref{balance})-(\ref{ohm-eq}) is invariant
under scaling of the $X$ vector; we scaled all reactances so that the largest
value was $100$ (also for solver robustness).  Interprocess communication in
our algorithm uses Unix sockets.  The computations below were performed on
three eight-core i7 machines with 48GB of RAM each.

\subsection{Second set of experiments}\label{experiments}

As before we use the Eastern Interconnect snapshot for our experiments. 
Synthetic contingencies were developed by removing $K$, random, highly
loaded lines, as in Section \ref{firstsetmethod}.

\noindent Our first set of experiments, shown in Table 1, concern cascades
with $R = 4$ rounds.   When applying rule (F.1), we used $\alpha = 0.55$.
As stated above, the segmented search was performed using $H = 50$ segments.

\begin{table}[h]
\centering
\caption{{\bf \emph{Performance of algorithm on 4-round cascades}}} 
\vskip 10pt
\begin{tabular}{|c | r| r| r|}
\hline
{\bf K} & {\bf yield,} &  {\bf yield,}  & {\bf wallclock} \\
 & {\bf no control} & {\bf control}  & {\bf (sec)} \\
\hline
\hline
{\bf 1} & 90.04  & 95.03 & 268\\
{\bf 2} & 1.25  & 50.13 & 174\\
{\bf 5} & 32.94 & 81.05 & 214\\
{\bf 10} & 2.02  & 36.97 & 194\\
{\bf 20} & 1.64  & 27.84 & 220\\
{\bf 50} & 0.83  & 16.96 & 477 \\
\hline
\end{tabular}
\label{table1}
\end{table}

\noindent In this table, the columns headed 'yield' indicate the 
{\em percentage} of total initial demand satisfied at the end of the cascade (without control, and using 
the computed control), and 'wallclock' is the observed 
parallel running time of the method.  In each of these runs, the total number
of gradient steps was small, typically smaller than 5.  

Note that in the
case $K = 1$ the interdiction has limited effect, but even so the control
is able to recover additional demand.  In the case $K = 5$ the
demand loss in the no-control case is substantial, but so is the benefit
of the control.  Finally, in the cases $K = 2, 10, 20, 50$ the network collapses
but the control sill recovers a significant amount of demand.  More experiments of this type will be forthcoming.

In the next set of experiments we use the case $K = 50$ in Table \ref{table1} 
to investigate in more detail 
the behavior of the algorithm as $R$ increases.  We used $\alpha = 0.5$ for
all these experiments.  Note that keeping $\alpha$ constant but increasing
$R$ effectively considers cascades that take longer from a 'real time' 
perspective, thereby giving more power to an agent applying control. 
If, instead, we were to increase $R$ while also decreasing $\alpha$, thus
giving more weight to 'history', we would be able to model cascades that
last for a fixed period of time, but where the individual rounds encompass
shorter spans of time.

Table \ref{tablerounds} reports on the experiments. As before, 'yield' 
is the percentage of demand satisfied at the end of the cascade, using 
no control, the control obtained by grid-search, and the control obtained
by segmented gradient search (started at the control computed by grid search).
The two wallclock columns report, in seconds, the time used by grid- and 
gradient-search.  In the case of grid search, we report the time spent on each
of the two search steps (i.e., over rounds 1 and 2, respectively). 
The column labeled 'grad steps' reports the number of 
gradient steps.

\begin{table}[h]
\centering
\caption{{\bf \emph{Impact of increasing number of rounds on K = 50 case from Table \ref{table1}}}} 
\vskip 10pt
\begin{tabular}{|c | r| r | r| r| r| r|}
\hline
{\bf R} & {\bf yield} &  {\bf yield} & {\bf yield}  & {\bf wallclock} & {\bf wallclock} & {\bf grad} \\
 & {\bf no control} & {\bf grid}  & {\bf gradient}& {\bf grid}  & {\bf gradient} & {\bf steps} \\
\hline
\hline
{\bf 5} & 4.13  & 18.11 & 31.86 & 30 + 17 & 1340 & 7\\
{\bf 6} & 2.02 & 23.01 & 25.86 & 26 + 14 & 657 & 6\\
{\bf 7} & 2.25  & 25.10 & 25.98 & 33 + 15 & 434  & 3\\
{\bf 8} & 0.78  & 29.27 & 46.97 & 18 + 43 & 3151 & 10 \\
\hline
\end{tabular}
\label{tablerounds}
\end{table}

\noindent Next we will comment on Table \ref{tablerounds}.\\

\noindent {\bf Computational workload} 

\noindent Consider the
case $R = 8$. Since we are using 
$H = 50$ segments, we have altogether $100$ control variables $c^r_i$ and $s^r_i$ per round $r$.
Since there are $7$ rounds during which we will apply control, we have a total of
$700$ individual variables.  Each partial derivative estimation requires
two simulations; thus in total each gradient estimation entails
$1400$ cascade simulations.   Per iteration, the
step-size computations  require $200$ additional cascade simulations; for 
a total  of $1600$ simulations per iteration of Procedure \ref{stochgrad1}.
The case $R = 8$ required
$10$ gradient iterations, and thus in total $16000$ simulations.  Each 8-round
 simulation (of the $15000$-bus Eastern
Interconnect, and using one core of the i7 CPU)
requires, on average, $4.5$ CPU seconds. This is primarily due to the
two power flow computations
per round, and linear solver data structure cleanup at the end of 
the simulation (and to a much lesser degree, to graph algorithms used to 
identify islands). 
Thus in total the computation
of the $R = 8$ case required approximately $72000$ CPU seconds. 
Since we have 24 worker cores, this 
translates to approximately $3000$ wallclock seconds. The balance of time
with respect to the actual wallclock time in Table \ref{tablerounds} (i.e., $151$ seconds) is due to inter-process
communication and networking delays, and logging of statistics to disc by the master.  On a per-simulation, per-core basis, this amounts
to $151*24/16000 \approx 0.22$ seconds, or roughly $5 \%$ as compared to
$4.5$ seconds total per simulation.  \\

\noindent{\bf Grid-search vs gradient-search}

\noindent In several cases gradient search significantly improves on the grid
search solution.  This is especially noticeable in the $R = 8$ case, and we
will examine this case in some detail.  

First, the grid-search control we
computed in this case uses $(\bar c^1, \bar b^1, \bar s^1) = (1, 1, 0.00018)$,
and $( \bar c^r, \bar b^r, \bar s^r) = (1,1,0)$ for $r > 2$, effectively 
limiting control to the first round.  In contrast, the gradient-search control
we computed applies control as late as round 7
(which, as we have 8 rounds in total, is the last round for which we 
compute a control as per our control template (\ref{affinecontcasc})) , and within earlier rounds
it applies different controls to different segments. In particular, in round
1 the gradient-search control uses the control vector $(0.95, 1, -0.0499)$
for segment 1 (the highest demand segment) as well as two other segments,
while for all other segments it uses $(1, 1, 0.00018)$.  And in round 7
it uses $(1.1, 1, 0.05)$ for
segment 2, and $(0.95, 1, 0.05)$ for segment 3, while for all other segments
it uses $(1, 1, 0)$.  Other controls different from $(1, 1, 0)$ are used
in rounds 2 and 4, while on rounds 5 and 6 it uses (1,1,0) throughout.

\begin{table}[h]
\centering
\caption{{\bf \emph{Controlled cascade evolutions}}} 
\vskip 10 pt
\begin{tabular}{|c| r r r r| r r r r|}
\hline
\multicolumn{1}{|c|} {} & \multicolumn{4}{c|}{{\bf Grid-search}} & \multicolumn{4}{c|}{{\bf Gradient-search}}\\
{\bf Round} & \boldmath{$\kappa$} & {\bf faults} & {\bf comps} & {\bf yield}  & \boldmath{$\kappa$} & {\bf faults} & {\bf comps} & {\bf yield} \\\hline \hline
{\bf 1} & 3.79 & 126 & 1 & 45.37 & 172.22 & 1629 & 32 & 60.72 \\  
{\bf 2} & 33.49 & 32 & 1 & 45.37 &  97.44 & 1079 & 293 & 54.26 \\  
{\bf 3} & 7.44  & 26 & 2 & 45.27 & 59.97 & 282 & 401 & 49.87 \\  
{\bf 4} & 6.69  & 82 & 4 & 45.27 & 21.88 & 89 & 459 & 48.67 \\  
{\bf 5} & 4.95  & 72 & 9 & 45.23 & 2.74 & 55 & 468 & 47.74 \\  
{\bf 6} & 1.99  & 28 & 13 & 45.23 & 13.27 & 10 & 471 & 47.72 \\  
{\bf 7} & 1.54  & 16 & 13 & 45.23 & 1.01 & 14 & 478 & 47.41 \\  
{\bf 8} & 1.00  & 16 & 13 & 29.27 & 1.00 & 1 & 478 & 46.97 \\  
\hline
\end{tabular}
\label{gridgradcomparo}
\end{table}

Table \ref{gridgradcomparo} compares the cascade evolution under both controls.  
Here, we report the value, at the end of each round of: $\kappa$ (the maximum line overload); the number of components of the network; and the yield.  'faults'
is the number of line outages experienced during each round.  Comparing the two evolutions, we see that the gradient-search control allows significantly more
outages in initial rounds (as well as much more islanding).  Nevertheless, the
number of outages is (nearly) monotonically decreasing under gradient-search and
by round 5 it is smaller than under grid-search.  Thus the gradient search
control appears to be maintaining high yield while carefully allowing outages 
to take place.  During round 7, gradient-search makes a large improvement on
the maximum line overload (from $13.27$ to $1.01$) but {\em without} losing
much yield; this is evidence of yet more sacrificial line outages.\\

\noindent {\bf Qualitative issues}

\noindent A comparison of 
the entry in Table \ref{tablerounds} for $R = 5$ and those for $R = 6, 7$
might appear to indicate that the controls computed for $R = 6$ and  $7$
are locally optimal,
because the 
control
that achieves yield $31.86\%$ for $R = 5$ ``should be'' 
feasible for all $R \ge 5$.

While it is true that the controls in Table \ref{tablerounds} 
can all be improved upon, the argument in the above paragraph is not
quite correct.  Refer to our Cascade Control template \ref{contcasc}.  The termination
step constitutes a last-recourse form of control -- if there are line
overloads at the start of the last round, loads are scaled so as to remove
the overloads, and in that case the cascade is considered terminated, regardless
of history (and of particular, of rule (\ref{memory-1}).  We model 
termination this way on purpose, so as to provide an agnostic termination
criterion that does not depend on numerical parameters of our model, in
particular, $\alpha$. Consider Table \ref{manykinds}.

\begin{table}[h]
\centering
\caption{{\bf \emph{Maximum line overload at end of each round for K = 50 case from Table \ref{tablerounds}}}} 
\vskip 2pt
\begin{tabular}{|c |  r r r r r|}
\hline
{\bf } & \boldmath$C_5$ & \boldmath$C_6$  &  \boldmath$C_7$ & \boldmath$C_8$ & {\bf None} \\ 
\hline
\noalign{\smallskip}
{\bf 1} & 6.47  & 1.83 & 2.22 & 3.79 & 177.83 \\
{\bf 2} & 14.12 & 1.83 & 1.57 & 33.49 & 122.06 \\
{\bf 3} & 36.79 & 1.23 & 1.30 & 6.90 & 114.45 \\
{\bf 4} & 1.72  & 1.14 & 2.26 & 6.70 & 22.47 \\
{\bf 5} &       & 0.99 & 1.18 & 59.33 & 45.43 \\
{\bf 6} &       &      & 1.08 & 1.98 & 40.33 \\
{\bf 7} &       &      &      & 1.18 & 114.90 \\
\hline
\end{tabular}
\label{manykinds}
\end{table}
In this table, the columns labeled ``$C_k$'', for $k = 5, \ldots, 8$  represent
the controls in Table \ref{tablerounds}, whereas ``None'' means no control.
The table shows, for each round, the maximum line overload at the {\em end} of 
that round, for each control option.  We see that $C_5$ reaches the start of
the termination round, round $5$, with maximum overload $1.7232$; the 
current yield at the start of round $5$ is $54.90 \%$ (not shown in the Table)
and most of the demand is in one island.  Hence the termination step will 
scale demands by $1/1.7232$ and yield will drop to $54.90/1.7232 \approx 31.86$,
 as Table \ref{tablerounds} shows.  As per our rules, this terminates
the cascade, although since $\alpha = 0.5$, 
and because the end-or-round $3$ maximum overload is very large, the 
maximum history-dependent line overload will be much larger than $1.732$
(it should be at least $0.5*36.79 = 18.40$.  Hence control $C_5$, if implemented
in a cascade with $6$ or more rounds, will not result in a stable state by
the end of round $5$.

Another point that emerges from Table \ref{manykinds} is that 
$C_5$ and $C_8$ 
tend to maintain higher line overloads late into the cascade --  this 
is a severe cascade, and having more time to apply control pays off.
But by doing so $C_5$ and $C_8$  are likely less robust.  Rather than performing
the same robustness analysis as in Section \ref{firstset}, we will next
consider the stability of the controls with respect to the $\alpha$ 
parameter in eq. (\ref{memory-1}) which in the above tests was set to
$0.5$.

This is a delicate issue, because the value of $\alpha$ is related to the
time duration of a round, and thus the structure of an optimal control 
{\em should}
depend on $\alpha$ (in other words, how much time we have impacts what 
kind of control we can apply).  The question is how stable a control remains as
$\alpha$ is perturbed.

\begin{table}[h]
\centering
\caption{{\bf \emph{Stability of controls in Table \ref{tablerounds} as a function of $\alpha$}}} 
\vskip 2pt
\begin{tabular}{|c |  r r r r |}
\hline
\boldmath$\alpha$ & \boldmath$C_5$ & \boldmath$C_6$  &  \boldmath$C_7$ & \boldmath$C_8$  \\ 
\hline
\noalign{\smallskip}
{\bf 0.45} & 1.49  &  25.05 &24.52  & 27.10  \\
{\bf 0.46} & 1.49  & 25.33   &24.52  &25.31   \\
{\bf 0.47} & 28.49 & 25.33 &24.52  & 25.31  \\
{\bf 0.48} & 28.47  &25.33  &24.52  & 26.08  \\
{\bf 0.49} & 28.47  & 25.33 &24.52  & 28.56  \\
{\bf 0.50} & 31.86  & 25.86 & 25.98 & 37.72  \\
{\bf 0.51} & 21.99  & 25.86 &25.98  & 34.11  \\
{\bf 0.52} & 20.99  & 25.86 &25.98  &  35.94 \\
{\bf 0.53} & 20.99  &25.86  &25.98  & 32.75  \\
{\bf 0.54} & 20.99  & 25.86 &25.98  & 32.75  \\
{\bf 0.55} & 20.99  & 25.86 & 25.98 & 31.83   \\
\hline
\end{tabular}
\label{alphadep}
\end{table}

In Table \ref{alphadep} we show the yields obtained by
running the $C_k$ controls from Table \ref{tablerounds} using their
respective numbers of rounds, but using different values for 
$\alpha$.  We see that in terms of the deviation from the nominal
case (i.e., $\alpha = 0.5$), $C_6$ and $C_7$ prove the most stable,
$C_8$ significantly less so and $C_6$ is very unstable.  It is
still the case that $C_8$ remains best overall: this is due to 
the severity of the cascade.

\begin{figure}[htb] 
\centering
\begin{center}
\includegraphics[height=4in, angle=270]{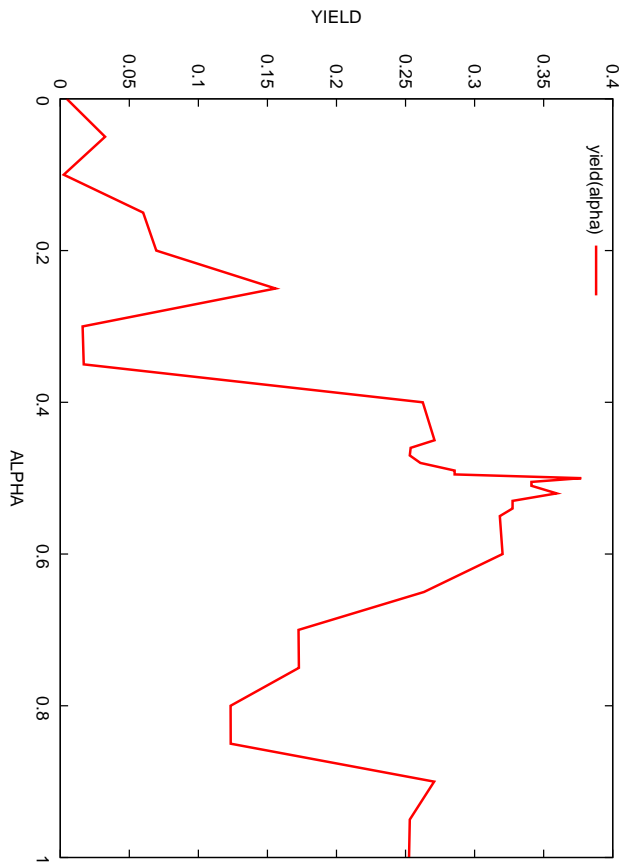}
\caption{yield as a function of alpha}\label{alphavsyield}
\end{center}
\end{figure}

In Figure \ref{alphavsyield} we consider a broader experiment along the
same lines.  For this experiment we used a control intended for the 
$R = 8$, $K = 50$ case considered in previous sections, and computed 
assuming $\alpha = 0.5$ (as per the memory-dependent rule (\ref{memory-1}) for
arc outages).  

Figure \ref{alphavsyield} displays the actual yield as $\alpha$ is changed 
away from $0.5$.  We note that yield is relatively robust for $\alpha$ larger
than but close to $0.5$, but not if $\alpha$ is decreased.  
We believe that this behavior is due to application of our (deterministic)
control results on lines becoming $100\%$ loaded.  Several heuristics,
based on ``padding'' (increasing) or ``tightening'' (reducing)
the flow limits $u_j$, suggest themselves.  However the greedy nature of
deterministic controls will remain a fundamental difficulty.  Section
\ref{stochasticoptimization} discusses the computation of controls 
under stochastics.

\subsubsection{Additional tests}\label{additional}
\noindent The next set of experiments address basic conjectures that arise
in the context of the type of control that we consider:
\begin{itemize}
\item It is best to stop the cascade in the first round, i.e. to sufficiently
reduce demands in the first round so as to eliminate all line overloads.
\item It is best to apply control in the first round only, and ride out the cascade for the remaining rounds.
\end{itemize}
\noindent In fact it is a simple task to create small examples where both
conjectures above are proved wrong.  Instead we explore these questions 
using the Eastern Interconnect, with a random interdiction 
of the type described above with $K = 50$, three rounds, and $\alpha = 1$ (no memory, and thus we obtain outage rule (F.1)). The results of this experiments are
as follows:
\begin{itemize}
\item Using no control, after three rounds $0.63 \%$ of the demand is satisfied.
\item Grid search produces a control with yield $45.46 \%$.
\item Starting from this control, and using segmented search with $50$ segments
improves yield to $50.02 \%$.
\end{itemize}
To gain a different perspective on this cascade, consider Table \ref{tableover},
which shows the distribution of line overloads at the start of the cascade.

\begin{table}[h]
\centering
\caption{{\bf \emph{Distribution of line overloads}}} 
\vskip 10pt
\setlength{\extrarowheight}{3pt}
\begin{tabular}{|c || r | r| r| r| r| r| r| r| r| r| r| r| r|}
\hline
{\bf $\lceil$ overload $\rceil$} & 1505 & 58 & 48  & 32 & 22 & 19 & 11 & 7 & 6 & 5 & 4 & 3 & 2 \\
{\bf count} & 1 & 1 & 2  & 1 & 2 & 1 & 1 & 2 & 2 & 4 & 6 & 18 & 181 \\
\hline
\end{tabular}
\label{tableover}
\end{table}

Where $f^1$ denotes the flow vector at the start of the cascade, the table
indicates the quantity of lines $j$ whose (rounded-up) overload 
$\lceil |f^1_j|/u_j \rceil$ equals a particular value.  Thus, $181$ lines 
$j$ are such that $1 <  |f^1_j|/u_j  \le 2$, $18$ lines $j$ are such
that $2 <  |f^1_j|/u_j  \le 3$, and so on.  One line has overload greater than
$1504.93$.  We will provide a more detailed analysis of this case in the near future; 
however it is the case (as may seem plausible from the table) that in order
to stop {\em all} overloads in round 1 a drastic reduction of demands is
needed.  We will instead describe the behavior of the optimal control 
computed by grid search has
$(c^1_v, b^1_v, s^1_v) = (300.7, 1, 0.0004)$, and 
$(c^2_v, b^2_v, s^2_v) = (1.26, 1, 0.62)$.  

Thus, in round 1, all demands
will be scaled by a (approximately) factor of $1.0 + 0.0004*(300.7 - 1504.93) =
0.51831$.  Considering Table \ref{tableover}, we see that in round 1 all lines
included in the columns with overload greater than $2$ will be outaged --
this is a total of $41$ lines, and in fact three more with overload close to
$2$ are outaged.  

At the start of round 2, the maximum overload is approximately $1.36194$.  Thus,
the control specifies that demands will be reduced by a factor of
$1.0 + 0.62*(1.26 - 1.36194) = 0.9368$.  This does not completely remove
all overloads and an additional $4$ lines become outaged.  

Finally, at the start of round 3, the maximum overload is $1.067891$.  By the
rules of our cascade model, this overload is now removed by scaling down
all demands.  Altogether, we obtain a yield of $0.51831*0.9368/1.067891 = 0.4547$, as stated above.  

\section{Stochastic optimization}\label{stochasticoptimization}
To further motivate the need for stochastic modeling, we consider the
same setup as for the experiment in Figure \ref{alphavsyield}: 
$R = 8$, $K = 50$ and $\alpha = 0.5$.  We 
simulated the behavior of the computed control
using the stochastic outage rule (\ref{actual1})-(\ref{actual3}), repeated
here for convenience.  We are given a parameter $0 < \epsilon < 1$; if $\tilde f$ is a flow vector, then line $j$
is {\em not} outaged if $|\tilde f_j| < (1 - \epsilon)u_j$, it
{\em is} outaged if  if $|\tilde f_j| > u_j$, and is outaged with probability
$1/2$ if $(1 - \epsilon)u_j \le \tilde f_j \le u_j$.

\begin{figure}[htb] 
\centering
\begin{center}
\includegraphics[height=5in, angle=270]{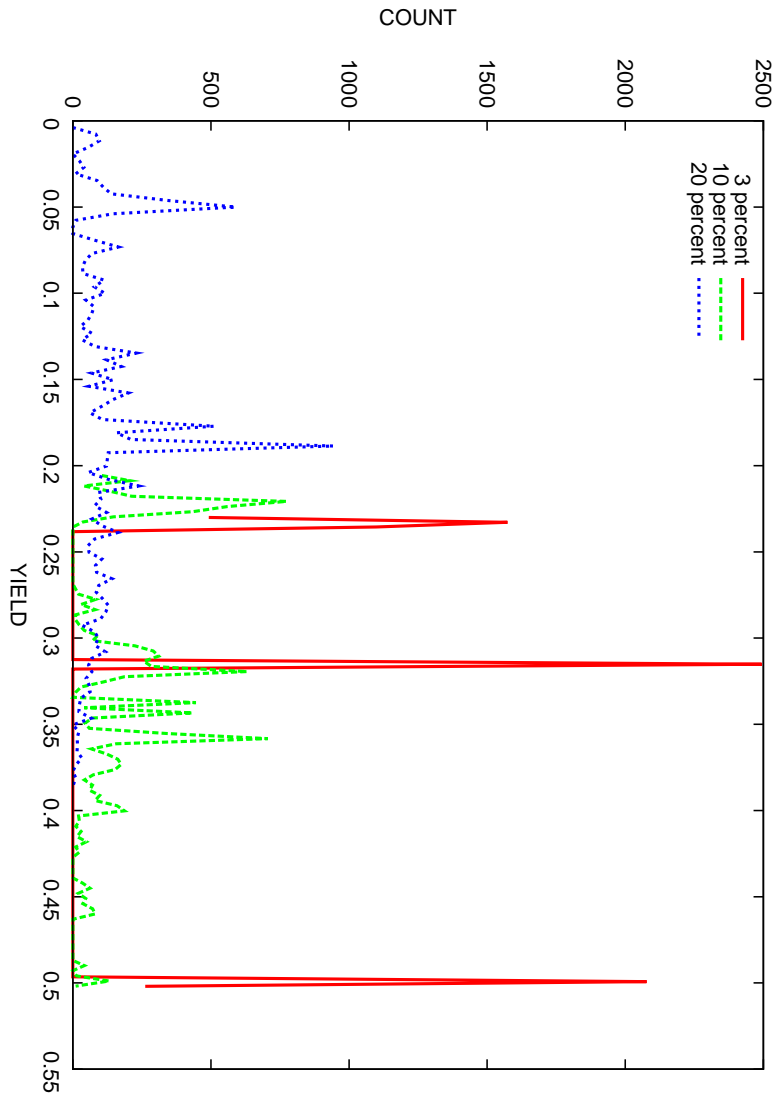}
\caption{yield histogram under stochastic outages}\label{epsilon}
\end{center}
\end{figure}

For various values of the tolerance $\epsilon$, we simulated the cascade
$10000$ times.  For $\epsilon < 0.03$ little difference was observed with
the nominal (deterministic) setting in that the average yield was close
to (the deterministic yield of) $50 \%$.  For $\epsilon = 0.03, 0.10$
and $0.20$ the results are displayed in Figure \ref{epsilon}.

The figure shows, for each value of the yield (and for each of the three
choices for $\epsilon$) how frequently that yield was observed.  Note that
for $\epsilon = 0.03$ the distribution is essentially trimodal.  This is typical
behavior and it points to a small number of critical lines,
which, in turn, result in a small number of cascade trajectories being overwhelmingly
likely.  For $\epsilon = 0.20$
yields close to zero are observed, but, significantly, the average is positive.

In a stochastic setting, the yield $\Theta^R(c,b,s)$ resulting from 
a control $(c,b,s)$ is a random variable.  
Below we discuss different methodologies for computing a (locally) 
optimal stochastic control, with the objective of maximizing 
the expectation $\rm E( \Theta^R)$.   Other merit criteria are also
of interest (and possibly better), such as a linear combination of expectation
and variance $\rm E(\Theta^R)  -  \lambda \, var(\Theta^R)$ 
($\lambda \ge 0$), or a Sharpe-ratio quantity $\rm E(\Theta^R)/var(\Theta^R)$.

The computational challenges inherent in any of these tasks are 
significant.  First, as displayed in Figure \ref{epsilon}, yield variances
can be extremely large.  From a theoretical perspective, the number of samples
needed to obtain reliable estimates become inordinately large.  Additionally,
there are subtle difficulties due to the non-monotonic behavior of power
flow systems (see \cite{bie2}), which is reminiscent of Braess' paradox \cite{braess}.

An example of this behavior is provide by our stochastic outage
rule (\ref{stochout}).  Note that under this rule, a line is more likely
to become outaged than under the deterministic rule (i.e., 
when $(1 - \epsilon_r)u_{j} < \tilde f^{r}_{j} \le u_{j} $ the line may
become outaged in the stochastic setting; not so in the deterministic 
setting).  Nevertheless, one can produce cases where the deterministic yield
of a control is {\em smaller} than a sample yield of the same control under
rule \ref{stochout}.  This phenomenon slows down convergence of our algorithms
and makes algorithm calibration difficult.

\subsection{Optimization methods through simulation}\label{gridstoch}

In either the first-order procedure (\ref{stochgrad1}) discussed above,
or in the grid-search setting, we obtain
a counterpart valid for a stochastic setting by replacing each (deterministic)
evaluation of a yield $\Theta^R(c,b,s)$ by an estimation of 
$\rm E ~ \Theta^R(c,b,s)$.  

This is the approach used in the next set of experiments, which parallel
those in Section \ref{firstset}.  For convenience, we restate the setup for
these tests.

First, two random, though adversarially chosen 
lines were removed from the grid.  Next, we compute the best control such
that
\begin{itemize}
\item [(i)] $c^r_v = b^r_v = 1$ for all $v$ and $r$.
\item [(ii)] $s^r_v = 0$ for all $v$ and $10 < r$. Thus, no control is applied
after round 10.
\item [(iii)] For each $1 \le r \le 10$, \underline{either} $s^r_v = 0.005$ for all
$v$, \underline{or} $s^r_v = 0$ for all $v$.
\end{itemize}
This was done, in the deterministic setting, while setting the maximum number of rounds, $R$, to $10$, $15$ 
$20$ and $25$.  In Section \ref{firstset}, we observed that the cascade is 
characterized by high line overloads during the initial rounds; nevertheless
if ``enough'' rounds are allowed ($34$) then without 
control the cascade stabilizes and produces approximately $78\%$ yield,
and this was slightly superior to what was obtained by the controls we 
computed.  In
summary, this case provides a good contrast between the (opposing) need to
``wait out'' an initially severe cascade on the one hand, with uncertainty
growth as we increase the number of rounds. Using the framework of stochastic outage rule (\ref{stochout}), 
with 
\begin{eqnarray}
\epsilon_r & = & 0.01 \, + 0.05* \lfloor r/10 \rfloor, \label{noise1prime}
\end{eqnarray}
we observed that the (deterministic) $20$-round control appeared superior.

In this section we will repeat these comparisons, except that now we will
compute controls of the form (i)-(iii) that (approximately) maximize
the expected yield. We compute such controls by modifying our
grid search: we evaluate a control vector by simulating $20$
cascades and computing the sample average yield.  This can be a somewhat
coarse approximation because, depending on the model for noise,
many more than $20$ samples may be needed for an accurate answer since the
variance of yield can be high.

\begin{table}[h]
\centering
\caption{{\bf \emph{Stochastic grid search results in case from Section \ref{firstset}}}} 
\vskip 10pt
\setlength{\extrarowheight}{3pt}
\begin{tabular}{|c || c | c| c| c|| c| c|}
\hline
{\bf R} & {\bf ave kR} & {\bf std kR} & {\bf ave kR}  & {\bf std kR} & {\bf ave cR} & {\bf std cR}\\
 & {\bf N = 20} & {\bf N = 20} & {\bf N = 1000}  & {\bf N = 1000} & {\bf N = 1000} & {\bf N = 1000} \\
\hline
{\bf 20} & 55.78 & 11.88 & 50.23  & 18.57 & 41.90 & 27.47 \\
{\bf 15} & 48.85 & 10.03 & 41.32  & 13.43 & 33.94 & 22.51 \\
{\bf 10} & 37.16 & 10.74 & 28.65  & 15.13 & 7.54 & 9.55 \\
\hline
\end{tabular}
\label{tablestochgrid}
\end{table}

For $R = 10, 15, 20$, denote by kR the control computed by the algorithm.
Table \ref{tablestochgrid} shows the results of our experiments. 
In the table, for each $R$,
'ave kR' and 'std kR' are the estimates of the average and standard 
of the yield of kR using $N = 20$ and $N = 1000$ samples.  Finally, 'ave cR' and 'std cR' are
the 1000-sample average and standard deviation of yield of the {\em deterministic} controls c20, c15, c10, computed in Section \ref{firstset} (as in Table \ref{table4}). 

We see that the kR controls are uniformly superior to their cR counterparts,
sometimes by almost one standard deviation.  We observe  that 
k20 is superior to k15, and much superior to k10.  This parallels observations
made in Section \ref{firstset}.

\subsection{Stochastic gradients}
The stochastic gradients method is a well-known approach for solving
optimization problems with uncertainty.  Because of the nonconcave nature
of the yield maximization problem, in our case it will amount to a local
search method.  Furthermore, there are some difficulties that are 
caused by the nonsmoothness in our models.
We will only outline here how we are approaching these difficulties.

The core step in the stochastic gradients approach is to (randomly)
sample a cascade, and, keeping the
cascade fixed, to compute the impact on yield of infinitesimally small
changes in the control parameters.  This is followed by a line search
to optimize the step size. This basic step is repeated.

A difficulty that we encounter when we attempt to 
make this outline more specific is that yield is not a differentiable 
function of the control parameters, for several reasons, the main one being the
stochastic outage protocol (\ref{stochout}), which, while smoother than a 
completely deterministic rule, is not smooth enough, due to its abrupt
transition between stochastic and deterministic regimes.  

We modify rule (\ref{stochout}) so that the probability of a  line 
will outage is always strictly positive and strictly smaller than $1$; however
when the overload is larger than $1$ the outaging probability will be very
close to $1$, and when the overload is significantly less than $1$ the
outaging probability (while positive) will be very small.  To this effect
consider a function 
$$ F \, : \, \field{R}_+ \, \rightarrow \, [0,1), \ \ \ \mbox{with} \ \ F(0) = 0 \ \mbox{and} \ F(x) \rightarrow 1 \ \ \mbox{as} \ \ x \rightarrow +\infty,$$
where the convergence is very rapid.  An example is $F(x) \, = \, 1 - e^{-Mx}$, for large $M > 0$. Likewise, consider a function
$$ G \, : \, [0,1] \, \rightarrow \, [0,1), \ \ \ \mbox{with} \ \ G(0) = 1 \ \mbox{and} \ G(x) \rightarrow 0 \ \ \mbox{as} \ \ x \rightarrow 1,$$
and again with rapid convergence. An example is  $G(x) \, = \, e^{-Mx}$ 
for large $M > 0$.

Having chosen $F$ and $G$, we modify outage rule (\ref{stochout}) as follows.
At round $r$, and given a tolerance $0 \ge \epsilon_r < 1$,
line $j$ is outaged
\begin{eqnarray}
\mbox{with probability} && \left \{
\begin{array}{lll}
\frac{1}{2}G\left(1 - \tilde f^r_{j}/(1 - \epsilon_r) \right), & \mbox{if $\tilde f^r_{j} \le (1 - \epsilon_r)u_{j}$} & \\
& & \\
\frac{1}{2}, & \mbox{if} \ {(1 - \epsilon_r)u_{j} < \tilde f^{r}_{j} < u_{j}} &  \label{smoothstoch} \\
& & \\
\frac{1}{2} \left(1 + F\left(\tilde f^{r}_{j}/u_{j} \ - \ 1 \right) \right),  & \mbox{if $u_{j} \le \tilde f^{r}_{j}$}. & \\
\end{array} \right. \label{hehehe}
\end{eqnarray}
In other words, if $\tilde f^{r}_j > u_j$, the outage probability 
is very large, but is bounded strictly away from $1$, and if
$\tilde f^{r}_j < (1 - \epsilon_r)u_j$ the outage probability is very small
but remains strictly positive.  By choosing $F$ and $G$ appropriately we
obtain an outage model that is arbitrarily close to rule (\ref{stochout}).

Another source of nonsmoothness in our models is the general form our 
control law in Step 2 of Procedure (\ref{affinecontcasc}).  However, it is
easy to see that the law can be approximated (arbitrarily closely) using
a smooth control.  

The computation of the (stochastic) gradient of the yield function at a 
given control vector $(\bar c, \bar b, \bar s)$ can now be described.  First, 
we sample a random cascade under the control $(\bar c, \bar b, \bar s)$ and 
outaging lines using rule (\ref{hehehe}).  This produces a particular 
sequence of lines that become outaged; i.e. at round $r$ a certain set
$S^r$ of lines is outaged, for $r = 1, 2, \ldots, R-1$.

Next, we compute the change in yield that results when we perturb the control
by a vector $(\epsilon^c, \epsilon^b, \epsilon^s)$ with infinitesimally
small entries, while still assuming
that set $S^r$ is the set of lines  outaged at round $r$, for each $r$.  This is a deterministic 
computation; rule (\ref{hehehe}) guarantees that the given cascade structure
retains positive probability.  This computation gives us the stochastic
gradient.

However, at this point we need to deal with the final reason that 
the yield function is not smooth, and this is the demand/supply 
adjustment in Step 3 of our generic cascade template (\ref{gencasc}) (or
in Step 6 of the cascade control template (\ref{contcasc})).  If, at round
$r$, under control $(c,b,s)$ a certain component $K$ has equal demand and
supply, then no adjustment takes place. 
However, even a small change in
the control that results in shedding less demand by round $r$ will not 
result in an increase of yield, whereas the opposite change in control will
possibly result in a decrease in yield.  Thus, in essence, a left derivative
is different from the right derivative; and moreover this is not a probability
zero event.

Nevertheless, it is still possible to adjust the stochastic gradients framework
so as to recover a valid first-order method.  The resulting approach is
related to the classical Frank-Wolfe method \cite{berts}.  In forthcoming work we will report on experiments with this approach.
\section{Upcoming work}
Our forthcoming work will focus on three areas: stochastics or robustness, in 
particular concerning the optimal scaling problem in Section \ref{scalingproblem}, an investigation of game-theoretic aspects of the type of control 
we study, and the use of AC power flow models.

With regards to the last point, a recent paper of Lavaei and Low 
\cite{lavaeilow} may yield a robust solver for AC power flow systems under
severe contingencies.  Even though the work in \cite{lavaeilow} 
relies on semi-definite programming, in fact one of the algorithms can
be restated as a second-order conic program, which may be more efficient.

\section*{Acknowledgment}
We would like to thank Ian Dobson and Ian Hiskens for fruitful discussions,
and for making the Eastern Interconnect data available to us.

\appendix
\section{Appendix - Formulations for the Optimal Demand Shedding Problem} \label{traditional}

We will first describe mixed-integer programming formulation for 
computing an optimal schedule for demand shedding, 
under deterministic line outages.  In the formulation below our
variables have the following interpretations:
\begin{itemize}
\item $f^r_{j}$ is the flow on line $(j)$ during round $r$ and $\pi^{r}_{j}$ 
($\nu^{r}_{j}$) is the positive (resp., negative) part of $f^{r}_{j}$
\item $\phi^r_{i}$ is the phase of bus $i$ during round $r$
\item for a demand bus $i$, $d^r_i$ indicates its demand during round $r$
\item for a generator bus $i$, $s^r_i$ indicates its supply during round $r$
\item for a line $j$, the variable $y^r_{j}$ takes value $1$ if
arc $j$ becomes outaged during round $r$ (and it takes value $0$ otherwise).
\item for a line $j$, the $0/1$ variable $p^r_{j}$ takes value $1$ if $f^r_{j} > 0$; likewise $n^r_{j}$ takes value $1$ if $f^r_{j} < 0$
\end{itemize}
For a demand bus $i$, we indicate by the constant $\tilde d_i$ its demand at the
start of the cascade.  Let $\tilde D$ denote the sum of all such quantities $\tilde d_i$.  
The formulation is as follows:
\newpage
\begin{eqnarray}
&& \max \sum_{i \in \cD} d^R_i \nonumber \\
\mbox{Subject to:} & & \sum_{j \in \delta^{+}(i)} f^r_{j} - \sum_{j \in \delta^{-}(i)} f^r_{j}  = \left \{
\begin{array}{lll}
~~s^r_i & {i \in \cG} &\\
-d^r_i &  {i \in \cD} &  \\
~~0   & \mbox{otherwise} &
\end{array} \right.  \ \forall \, 1 \le r \le R \label{1b}\\
&& \nonumber \\
& & f^r_{j} \, = \, \pi^r_{j} - \nu^r_{j} \,\,\, \forall \,\, j \in \cA \,\, \mbox{and $1 \le r \le R$} \label{defflow}\\
&& \nonumber \\
&& \pi^r_{j} \, \le \, \tilde D p^r_{j}, \ \  \nu^r_{j} \, \le \, \tilde D n^r_{j}, \,\,\, \forall \,\, j \in \cA \,\, \mbox{and $1 \le r \le R$} \label{defpq}\\
&& \nonumber \\
&&  p^r_{j} \, + \, n^r_{j} \, = \, 1 - \sum_{h=1}^{r-1} y^h_{j}, \,\,\, \forall \,\, j \in \cA \,\, \mbox{and $1 \le r \le R$} \label{sumpq}\\
&& \nonumber \\
& & \pi^r_{j} + \nu^r_{j} - u_{j} \, \le \, \tilde D y^r_{j} \,\,\, \forall \,\, j \in \cA \,\, \mbox{and $1 \le r \le R$} \label{defyr}\\
&& \nonumber \\
& & \pi^r_{j} + \nu^r_{j} \, \ge \, u_{j} y^r_{j} \,\,\, \forall \,\, j  \in \cA\,\, \mbox{and $1 \le r \le R-1$} \label{defyr,2}\\
&& \nonumber \\
& & \pi^R_{j} + \nu^R_{j}  \, \le \, u_{j} \,\,\, \forall \,\, j  \in \cA \label{finalround}\\
&& \nonumber \\
& & |\phi^r_i - \phi^r_j - x_{j}f^r_{j}| \, \le \, M_{j}~ \sum_{h=1}^{r-1} y^h_{j}   \ \ \ \forall j  \in \cA \label{ohm-eq-r} \\
&& \nonumber \\
& & 0 \le s_i^r \le \tilde s_i \,\, \, \forall \, i \in \cG, \ \ \ \  0 \le d_i^r \le \tilde d_i \, \, \, \forall \, i \in \cD,\\
&& \nonumber \\
&& p^r_{j}, \, n^r_{ij}, \, y^r_{j} \, = \, 0 \ \mbox{or $1$}, \,\,\, \forall \,\, j  \in \cA \,\, \mbox{and $1 \le r \le R$} \\
&& \nonumber \\
&& 0 \le \pi^r_{j}, \, 0\le \nu^r_{j}, \,\,\, \forall \,\, j  \in \cA\,\, \mbox{and $1 \le r \le R$}. 
\end{eqnarray}
In this formulation, (\ref{1b}) is a flow balance constraint: it specifies that
during round $r$ each generator $i$ outputs $s^r_i$ units of flow, and 
similarly with demand buses.  
Constraints (\ref{defflow})-(\ref{sumpq}) together with the fact that 
$p^r_{j}$ and $n^r_{j}$ are $0/1$ variables, guarantee that 
$\pi^r_{j}$ ($\nu^r_{j}$) is the positive (resp., negative) part of
$f^r_{j}$, and that furthermore $f^r_j = 0$ if line $j$ was outaged at a round
prior to $r$.  Constraint (\ref{defyr}) guarantees that if $y^r_{j} = 0$ then
$|f^r_{j}| \le u_{j}$ and (\ref{defyr,2}) guarantees that if 
$y^r_{j} = 1$ then $|f^r_{j}| \ge u_{j}$.    Thus, we obtain a
mix of the two alternative versions of rule (F.1). 

Constraint (\ref{finalround}) indicates the desired termination condition in
round $R$.  We note that (\ref{ohm-eq-r}) involves an absolute
value but is easily replaced by two standard linear inequalities.  The 
quantity $M_{j}$ is assumed to be a ``large enough'' quantity (see \cite{bie} 
for a related discussion).

Other formulations are possible, and in particular it is easy to enforce
additional rules constraining how demand can be shed.  
The formulation can also be adapted to use 
the memory-dependent outage models (\ref{memory-1}) or (\ref{memory-2}).
By adapting constraints (\ref{defyr}) and (\ref{defyr,2}) one can model
the stochastic rule (F.2), obtaining a (mixed-integer) stochastic program;
the underlying uncertainty is primarily of an endogenous nature; see
\cite{goelgrossman}.

\subsubsection{Discussion}

The above problem, even in its simplest, deterministic form is likely 
quite nontrivial.  Constraints (\ref{defyr}) and (\ref{defyr,2}) are essential
in modeling the outaging of lines; similar constraints are used in
formulations for the classical $N-K$ problem (see \cite{bie}, \cite{bie2}) 
and contribute to make the problem very difficult.  Note that we
would need to handle cases with $n > 10^4$, $m > 2\times 10^4$ (and $R > 2$).

A larger concern involves the fact that the optimal solution is likely to
entail very complex control strategies that may be difficult to implement.  
All the
approaches we discussed above, including the
stochastic programming versions of the formulation, are likely to 
specify very precise actions in each round (and scenario) which may be 
problematic in practice in what likely would be a very ``noisy'' environment.

Nevertheless, the study of the formulation may still prove a very worthwhile
exercise, one that could highlight underlying weaknesses of the models and
hidden vulnerabilities in the operation of the grid.

\end{document}